\documentclass[11pt, reqno, psamsfonts]{amsart}
\usepackage{anton}
\usepackage{jing}

\author{Anton~Bernshteyn}
\address{\normalfont (AB) Department of Mathematics, University of California, Los Angeles, CA, USA}
\email{bernshteyn@math.ucla.edu}

\author{Jing~Yu}
\address{\normalfont (JY) Shanghai Center for Mathematical Sciences, Fudan University, Shanghai, China}
\email{jyu@fudan.edu.cn}

\thanks{AB's research is partially supported by the NSF grant DMS-2045412 and the NSF CAREER grant DMS-2528522. JY's research was partially supported by ARC-ACO fellowship for the Spring 2022 semester.}
\title{\sffamily Large-scale geometry of Borel graphs of polynomial growth}

\begin{document}
    
    \maketitle

    
    \begin{abstract}
        We study graphs of polynomial growth from the perspective of asymptotic geometry and descriptive set theory. The starting point of our investigation is a theorem of Krauthgamer and Lee who showed that every connected graph of polynomial growth admits an injective contraction mapping to $(\Z^n, \|\cdot\|_\infty)$ for some $n \in \N$. We strengthen and generalize this result in a number of ways. In particular, answering a question of Papasoglu, we construct coarse embeddings  from graphs of polynomial growth to $\Z^n$. Moreover, we only require $n$ to be linear in the asymptotic polynomial growth rate of the graph; this confirms a conjecture of Levin and Linial, London, and Rabinovich ``in the asymptotic sense.'' (The exact form of the conjecture was refuted by Krauthgamer and Lee.) All our results are proved for Borel graphs, which allows us to settle a number of problems in descriptive combinatorics. Roughly, we prove that graphs generated by free Borel actions of $\Z^n$ are universal for the class of Borel graphs of polynomial growth. This provides a general method for extending results about $\Z^n$-actions to all Borel graphs of polynomial growth. For example, an immediate consequence of our main result is that all Borel graphs of polynomial growth are hyperfinite, which answers a well-known question in the area. As another illustration, we show that Borel graphs of polynomial growth support a certain combinatorial structure called a toast. An important technical tool in our arguments is the notion of padded decomposition from computer science, which is closely related to the concept of asymptotic dimension due to Gromov. Along the way we find an alternative, probabilistic proof of a theorem of Papasoglu that graphs of asymptotic polynomial growth rate $\rho < \infty$ have asymptotic dimension at most $\rho$ and establish the same bound in the Borel setting.
    \end{abstract}

    \tableofcontents
    \vspace*{-1cm}
        
    \section{Introduction and main results}

    
    \subsection{Graphs of polynomial growth} 
    
    This paper is a contribution to the study of large-scale geometry of graphs of polynomial growth. All graphs in this paper are undirected and simple, but they may \ep{and typically will} be infinite. We use standard asymptotic notation, in particular $O(\cdot)$, $\Omega(\cdot)$, and $\Theta(\cdot)$. Subscripts in the asymptotic notation, such as $O_{a,b,\ldots}(\cdot)$, indicate that the implied constants may depend on $a$, $b$, \ldots\ Throughout, we use the word ``metric'' to mean ``extended metric,'' i.e., we allow distances in a metric space to be infinite. Given a graph $G$, we let $\dist_G$ denote the \emphd{graph metric} on $V(G)$, that is, for two vertices $u$, $v \in V(G)$, $\dist_G(u,v)$ is the minimum number of edges on a $uv$-path in $G$ if such a path exists and $\infty$ otherwise. We also let $B_G(u,r)$ denote the closed ball of radius $r$ around a vertex $u$ in the metric $\dist_G$. The following definition introduces our central objects of study:
    
    \begin{defn}[{Graphs of polynomial growth}]
        Let $G$ be a graph. The \emphd{growth function} of $G$ is the mapping $\gamma_G \colon [0, +\infty) \to [0, +\infty]$ defined as follows:
        \[
            \gamma_G(r) \,\defeq\, \sup_{u \in V(G)} |B_G(u, r)|.
        \]
        For $r \geq 1$, we define the quantity $\rho(G, r)$ by the formula
        \[
            \rho(G, r) \,\defeq\, \frac{\log \gamma_G(r)}{\log(r+1)}.
        \]
        That is, we have $\gamma_G(r) = (r+1)^{\rho(G, r)}$ when $\gamma_G(r) < \infty$. The \emphd{exact growth rate} of $G$ is given by
        \[
            \er(G) \,\defeq\, \sup_{r \geq 1} \rho(G, r).
        \]
        The \emphd{asymptotic growth rate} $\ar(G)$ of $G$ is given by
        \[
            \ar(G) \,\defeq\, \limsup_{r \to \infty} \rho(G, r).
        \]
        If $\ar(G) < \infty$, we say that $G$ is a \emphdef{graph of polynomial growth}.
    \end{defn}
    
    It is clear from the definitions that $\er(G) \geq \ar(G)$ for any graph $G$, and this inequality may be strict. For example, if $G$ is finite, then $\ar(G) = 0$, while $\er(G)$ can be arbitrarily large. Notice, however, that if $\ar(G)$ is finite (i.e., if $G$ is of polynomial growth), then $\er(G)$ must be finite as well. In general, large-scale properties of a graph are better reflected in its asymptotic growth rate, as it ignores ``local'' phenomena that may influence the exact growth rate. 
    
\begin{figure}[t]
			\centering
			\begin{tikzpicture}[scale=0.8]
			    \begin{scope}[scale=2]
				    \filldraw (0,0) circle (2pt);
				    \filldraw (0.5,0) circle (2pt);
				    \filldraw (1,0) circle (2pt);
				    \filldraw (1.5,0) circle (2pt);
				    \filldraw (0,0.5) circle (2pt);
				    \filldraw (0.5,0.5) circle (2pt);
				    \filldraw (1,0.5) circle (2pt);
				    \filldraw (1.5,0.5) circle (2pt);
				    \filldraw (0,1) circle (2pt);
				    \filldraw (0.5,1) circle (2pt);
				    \filldraw (1,1) circle (2pt);
				    \filldraw (1.5,1) circle (2pt);
				    \filldraw (0,1.5) circle (2pt);
				    \filldraw (0.5,1.5) circle (2pt);
				    \filldraw (1,1.5) circle (2pt);
				    \filldraw (1.5,1.5) circle (2pt);
				    
				    \draw (-0.2,0) -- (1.7,0) (-0.2,0.5) -- (1.7,0.5) (-0.2,1) -- (1.7,1) (-0.2,1.5) -- (1.7,1.5) (0,-0.2) -- (0,1.7) (0.5,-0.2) -- (0.5, 1.7) (1,-0.2) -- (1,1.7) (1.5,-0.2) -- (1.5, 1.7);
				\end{scope}
				
				\begin{scope}[scale=2,xshift=4cm]
				    \filldraw (0,0) circle (2pt);
				    \filldraw (0.5,0) circle (2pt);
				    \filldraw (1,0) circle (2pt);
				    \filldraw (1.5,0) circle (2pt);
				    \filldraw (0,0.5) circle (2pt);
				    \filldraw (0.5,0.5) circle (2pt);
				    \filldraw (1,0.5) circle (2pt);
				    \filldraw (1.5,0.5) circle (2pt);
				    \filldraw (0,1) circle (2pt);
				    \filldraw (0.5,1) circle (2pt);
				    \filldraw (1,1) circle (2pt);
				    \filldraw (1.5,1) circle (2pt);
				    \filldraw (0,1.5) circle (2pt);
				    \filldraw (0.5,1.5) circle (2pt);
				    \filldraw (1,1.5) circle (2pt);
				    \filldraw (1.5,1.5) circle (2pt);
				    
				    \draw (-0.2,0) -- (1.7,0) (-0.2,0.5) -- (1.7,0.5) (-0.2,1) -- (1.7,1) (-0.2,1.5) -- (1.7,1.5) (0,-0.2) -- (0,1.7) (0.5,-0.2) -- (0.5, 1.7) (1,-0.2) -- (1,1.7) (1.5,-0.2) -- (1.5, 1.7);
				    
				    \draw (-0.2,-0.2) -- (1.7,1.7) (-0.2, 0.3) -- (1.2,1.7) (-0.2, 0.8) -- (0.7,1.7) (-0.2, 1.3) -- (0.2, 1.7) (0.3,-0.2) -- (1.7,1.2) (0.8,-0.2) -- (1.7,0.7) (1.3,-0.2) -- (1.7,0.2);
				    
				    \draw (1.5+0.2,-0.2) -- (1.5-1.7,1.7) (1.5+0.2, 0.3) -- (1.5-1.2,1.7) (1.5+0.2, 0.8) -- (1.5-0.7,1.7) (1.5+0.2, 1.3) -- (1.5-0.2, 1.7) (1.5-0.3,-0.2) -- (1.5-1.7,1.2) (1.5-0.8,-0.2) -- (1.5-1.7,0.7) (1.5-1.3,-0.2) -- (1.5-1.7,0.2);
				\end{scope}
			\end{tikzpicture}
			\caption{Fragments of the graphs $\gr_2$ (left) and $\gr_{2,\infty}$ (right).\label{fig:grids}}
	\end{figure}
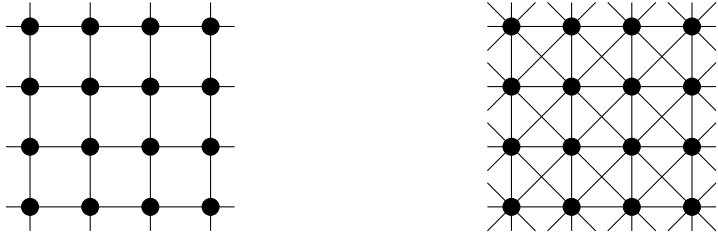

    Perhaps the most familiar example of a graph of polynomial growth is the \emphd{$n$-dimensional square grid} $\gr_{n}$, defined as follows. The vertex set of $\gr_{n}$ is the Abelian group $\Z^n$, and two vertices $u$, $v \in \Z^n$ are adjacent in $\gr_{n}$ if and only if $\|u - v\|_1 = 1$, i.e., if the difference between $u$ and $v$ is, up to sign, one of the standard generators for $\Z^n$. For our purposes, it will be more convenient to consider the closely related graph $\gr_{n,\infty}$, which we interchangeably call the \emphd{$n$-dimensional $\infty$-grid} or the \emphd{$n$-dimensional grid with diagonals}. The vertex set of $\gr_{n,\infty}$ is also $\Z^n$, and two vertices $u$ and $v$ are adjacent in $\gr_{n, \infty}$ if and only if $\|u - v \|_\infty = 1$ (see Fig.~\ref{fig:grids}). Note that both graphs $\gr_n$ and $\gr_{n,\infty}$ have asymptotic growth rate $n$ and exact growth rate $\Theta(n)$. For brevity, we write $\dist_1$ and $\dist_\infty$ to denote the graph metrics on $\gr_n$ and $\gr_{n,\infty}$ respectively.
    
    Levin and Linial, London, and Rabinovich conjectured that grid graphs (with diagonals) are, in some sense, universal among all connected graphs of polynomial growth \cite[Conjecture 8.2]{linial1995geometry}. To state their conjecture precisely, we need to introduce a few definitions. Recall that a \emphd{homomorphism} from a graph $G$ to a graph $H$ is a mapping $f \colon V(G) \to V(H)$ that sends edges to edges (i.e., if $uv \in E(G)$, then $f(u)f(v) \in E(H)$). The following definition is implicit in \cite{linial1995geometry}:
    
    \begin{defn}[{Injective dimension}]
        The \emphd{injective dimension} of a connected graph $G$, denoted by $\injdim(G)$, is the least $n \in \N$ such that $G$ admits an injective homomorphism to the graph $\gr_{n, \infty}$ if such $n$ exists; otherwise, $\injdim(G) \defeq \infty$.
    \end{defn}
    
    Clearly, if $G$ is a graph of injective dimension $n < \infty$, then $G$ must be of polynomial growth and, in fact, $\ar(G) \leq \ar(\gr_{n,\infty}) = n$ and $\er(G) \leq \er(\gr_{n,\infty}) = O(n)$. Levin and Linial, London, and Rabinovich conjectured that the converse to this statement also holds. More precisely, their hypothesis was that if $G$ is a connected graph with $\er(G) = \rho < \infty$, then:
    \begin{enumerate}[label=\ep{\normalfont{}H\arabic*}]
        \item $\injdim(G)$ is finite and bounded by a function of $\rho$;
        \item\label{item:H2} in fact, $\injdim(G) = O(\rho)$.
    \end{enumerate}
    Krauthgamer and Lee confirmed the first part of this conjecture and refuted the second:
    
    \begin{theo}[{Krauthgamer--Lee \cite[Theorems 5.8 and 3.2]{krauthgamer2007intrinsic}}]\label{theo:KL}\mbox{}
    
        \begin{enumerate}[label=\ep{\normalfont\arabic*}]
            \item\label{item:KL1} If $G$ is a connected graph\footnote{This result is stated for finite graphs in \cite{krauthgamer2007intrinsic}, but it generalizes to all connected graphs by compactness.} with $\er(G) = \rho < \infty$, then $\injdim(G) = O(\rho \log \rho)$.
            \item\label{item:KL2} Furthermore, the bound $O(\rho \log \rho)$ on the injective dimension is optimal. That is, for any $\rho \geq 1$, there exists a connected graph $G$ with $\er(G) \leq \rho$ and $\injdim(G) = \Omega(\rho \log \rho)$.
        \end{enumerate}
    \end{theo}
    
    \subsection{Informal summary of our main contributions}
    
    In this paper we strengthen and extend the first part of Theorem~\ref{theo:KL} in a number of ways:
    
    \begin{enumerate}[label=\ep{\emphd{\Alph*}}]
        \item\label{item:A} We prove an \emph{asymptotic} analog of Theorem~\ref{theo:KL}(1), which allows us to achieve the bound $O(\rho)$ on the dimension of the grid, where $\rho$ is the asymptotic growth rate of $G$.
        
        \item\label{item:B} We establish a more precise relationship between the graph metrics on $G$ and on the grid to which $G$ is mapped, giving an affirmative answer to a question of Papasoglu \cite{papasoglu2021polynomial}.
        
        \item\label{item:C} We extend Theorem~\ref{theo:KL} to the realm of Borel graphs and prove that graphs generated by free Borel actions of $\Z^n$ are, in an appropriate sense, universal for the class of Borel graphs of polynomial growth.
    \end{enumerate}
    
    \noindent Before stating our main results formally, let us briefly comment on each of these points.
    
    \subsubsection*{\ref{item:A} \normalfont{\itshape{}Using the asymptotic growth rate instead of the exact growth rate}}
    
    It is clear that the injective dimension of a graph constrains its exact growth rate, simply because if $\injdim(G) = n < \infty$, then $\gamma_G(r) \leq \gamma_{\gr_{n, \infty}}(r)$ for all $r$. Nevertheless, from the point of view of large-scale geometry, it seems more natural to treat the asymptotic growth rate as the relevant parameter. In order to eliminate the dependence on the local structure of $G$, we need to consider mappings that are not necessarily injective, but only injective in a suitable ``asymptotic'' sense. Namely, we shall allow vertices that are close to each other to be mapped to the same vertex. Formally, given a pair of graphs $G$, $H$, we say that a function $f \colon V(G) \to V(H)$ is:
    \begin{itemize}
        \item a \emphd{contraction} if $\dist_H(f(u), f(v)) \leq \dist_G(u,v)$ for all $u$, $v \in V(G)$, or, equivalently, if for all $uv \in E(G)$, we either have $f(u)=f(v)$ or $f(u)f(v) \in E(H)$;
        \item \emphd{asymptotically injective} if there is $r_0 \in \N$ such that $f(u) \neq f(v)$ whenever $\dist_G(u,v) > r_0$, or, equivalently, if the $f$-preimage of every vertex of $H$ is a set of diameter at most $r_0$ in $G$.
    \end{itemize}
    Now, given a graph $G$ of polynomial growth, we seek an asymptotically injective contraction from $G$ to the graph $\gr_{n, \infty}$, where $n$ is bounded by a function of $\ar(G)$. We prove that such a contraction exists and, moreover, one can take $n = O(\ar(G))$. This shows that part \ref{item:H2} of the Levin--Linial--London--Rabinovich conjecture, although false as stated, is ``true in the asymptotic sense.'' Note that this does not contradict Theorem~\ref{theo:KL}\ref{item:KL2}: the counterexamples constructed by Krauthgamer and Lee are finite graphs, and sending the entire vertex set of a finite connected graph to a single vertex trivially yields an asymptotically injective contraction.
    
    \subsubsection*{\ref{item:B} \normalfont{\itshape{}Relating the graph metrics on $G$ and on the grid}}
    
    If $f$ is an asymptotically injective contraction defined on a graph $G$, then for any pair of vertices $u$, $v \in V(G)$, the distance between $f(u)$ and $f(v)$ is bounded above by $\dist_G(u,v)$ and, moreover, it is bounded below by $1$ when $\dist_G(u,v)$ is large enough. At this point, it is natural to ask if better control over the distance between $f(u)$ and $f(v)$ is possible. In particular, can we bound it below by a function that tends to infinity with $\dist_G(u,v)$? This idea is captured by the notion of a coarse embedding, introduced by Gromov \cite{gromov1993asymptotic}:
    
    \begin{defn}[Coarse embeddings]\label{defn:coarse}
        Given a pair of metric spaces $(X, d_X)$, $(Y, d_Y)$, a mapping $f \colon X \to Y$ is called a \emphd{coarse embedding} if there exist non-decreasing functions $b_\ell$, $b_u \colon [0,\infty] \to [0,\infty]$ such that $b_\ell(r) \leq b_u(r) < \infty$ for all $r < \infty$, $\lim_{r \to \infty} b_\ell(r) = \lim_{r \to \infty} b_u(r) = \infty$, and 
    \[
        b_\ell(d_X(x,y)) \,\leq\, d_Y(f(x), f(y)) \,\leq\, b_u(d_X(x,y)) \qquad \text{for all $x$, $y \in X$}.
    \]
    \end{defn}
    Note that a coarse embedding from a graph $G$ to a graph $H$ must be asymptotically injective. Inspired by Theorem~\ref{theo:KL}, Papasoglu asked whether graphs of polynomial growth admit coarse embeddings into grid graphs \cite[Question 4.1]{papasoglu2021polynomial}. We answer this question in the affirmative. Furthermore, our coarse embeddings are contractions (i.e., one can take $b_u(r) = r$) and the function $b_\ell$ is almost linear. More precisely, we show that for any $\epsilon > 0$, a connected graph $G$ of polynomial growth admits a contraction $f$ to $\gr_{n, \infty}$ with $n = O_\epsilon(\ar(G))$ such that for all $u$, $v \in V(G)$,
    \begin{equation}\label{eq:LowerBound}
        \dist_\infty(f(u),f(v)) \,=\, \Omega_{G, \epsilon}(\dist_G(u,v)^{1-\epsilon}).
    \end{equation}
    Some dependence of the implied constants in \eqref{eq:LowerBound} on $G$ is unavoidable, since we do not know \emph{a priori} how large $r$ needs to be so that $\rho(G, r)$ can be bounded above in terms of $\ar(G)$. The version of this result that we actually prove is Theorem~\ref{theo:coarse_p}, which makes the bound completely explicit and is meaningful even when the graph $G$ is finite. Nevertheless, we have opted to use asymptotic formulas such as \eqref{eq:LowerBound} in the introduction for simplicity. 
    
    By analogy with Theorem~\ref{theo:KL}, for each $\epsilon > 0$, we also construct \emph{injective} coarse embeddings $f$ from a graph $G$ of polynomial growth into the $\infty$-grid $\gr_{n,\infty}$ with two different sets of parameters. In both cases, the lower bounding function $b_\ell$ satisfies $b_\ell(r) = \Omega_{G, \epsilon}(r^{1-\epsilon})$ \ep{that is, \eqref{eq:LowerBound} holds for $f$}, while the dimension $n$ of the grid and the upper bounding function $b_u$ are as follows:
    \begin{itemize}
        \item $b_u(r) = \max \set{r, O_{G, \epsilon}(1)}$ \ep{and hence $f$ is an \emphd{asymptotic contraction}, i.e., $\dist_\infty(f(u), f(v)) \leq \dist_G(u,v)$ whenever $\dist_G(u,v)$ is large enough} and $n = O_\epsilon(\ar(G))$;
        
        \item $b_u(r) = r$ (i.e., $f$ is a contraction) and $n = O_\epsilon(\er(G) \log (\er(G)))$.
    \end{itemize}
    The latter construction is a direct strengthening of Theorem~\ref{theo:KL}\ref{item:KL1}.

    Note that we cannot take $\epsilon = 0$ in the above results: Papasoglu constructed graphs of polynomial growth (indeed, of asymptotic growth rate arbitrarily close to $1$) whose so-called \emph{asymptotic Assouad--Nagata dimension} is infinite  \cite[Theorem~3.5]{papasoglu2021polynomial}, which in particular implies that they do not admit a quasi-isometric embedding into any finite-dimensional grid.
    
    \subsubsection*{\ref{item:C} \normalfont{\itshape{}Extension to Borel graphs}}
    
    Over the course of the past quarter century, a new and exciting research area has emerged at the intersection of combinatorics (especially graph theory) and descriptive set theory, referred to as \emphd{descriptive combinatorics}. For surveys of this subject, see \cite{KechrisMarks} by Kechris and Marks and \cite{pikhurko2020borel} by Pikhurko. Broadly speaking, descriptive combinatorics investigates infinite combinatorial objects endowed with an additional ``continuous'' structure, such as a topology, a measure, or a $\sigma$-algebra of Borel sets. The objects we are specifically interested in here are \emph{Borel graphs}:
    
    \begin{defn}[Borel graphs]\label{defn:BorelGraph}
        A \emphd{Borel graph} is a graph $G$ whose vertex set $V(G)$ is a standard Borel space and whose adjacency relation $\set{(u,v) \in V(G)^2 \,:\, uv \in E(G)}$ is a Borel subset of $V(G)^2$.
    \end{defn}
    
    We refer the reader to \cite{KechrisDST,AnushDST} for basic descriptive set-theoretic terminology, such as the definition of a standard Borel space. Note that every countable graph is trivially Borel, since every countable set in a standard Borel space is Borel. Hence, Definition~\ref{defn:BorelGraph} is only really interesting when $V(G)$ is uncountable. 
    By the Borel isomorphism theorem \cite[Theorem~15.6]{KechrisDST}, all uncountable standard Borel spaces are isomorphic to each other, so typically no generality is lost by assuming that $V(G)$ is some specific space, such as the unit interval $[0,1]$ or the Cantor space $2^\N$.
    
    In this paper we will be concerned with the case when $G$ has uncountably many connected components, while each individual component is countable. One context in which such graphs naturally arise is the study of finitely generated groups and their actions. Let $\G$ be a group generated by a finite symmetric\footnote{A subset of a group is \emphd{symmetric} if it is closed under taking inverses.} set $F \subseteq \G\setminus \set{\mathbf{1}_\G}$. The \emphd{Cayley graph} $\mathsf{Cay}(\G, F)$ of $\G$ corresponding to $F$ is the graph with vertex set $\G$ in which two vertices $\gamma$, $\delta \in \G$ are adjacent if and only if $\gamma \delta^{-1} \in F$. For instance, both $\gr_n$ and $\gr_{n,\infty}$ are Cayley graphs of the group $\Z^n$ corresponding to the generating sets $\set{z \in \Z^n \,:\, \|z\|_1 = 1}$ and $\set{z \in \Z^n \,:\, \|z\|_\infty = 1}$ respectively. 
    
    While $\mathsf{Cay}(\G, F)$ is a countable graph, there is a way to generalize its construction to obtain uncountable Borel graphs as well. To this end, let $\G \acts X$ be a Borel action of $\G$ on a standard Borel space $X$. Then we can form the so-called \emphd{Schreier graph} of this action by mapping the Cayley graph of $\G$ onto each orbit. More precisely, if $F \subseteq \G \setminus \set{\mathbf{1}_\G}$ is a finite symmetric generating set for $\G$, 
    then the corresponding Schreier graph $\mathsf{Sch}(X, F)$ has vertex set $X$ and edge set \[E(\mathsf{Sch}(X, F)) \,\defeq\, \set{\set{x, \sigma \cdot x} \,:\, x \in X,\ \sigma \in F,\ \sigma \cdot x \neq x}.\] \ep{The condition ``$\sigma \cdot x \neq x$'' is needed to ensure the graph has no loops}. It is clear from the definition that $\mathsf{Sch}(X, F)$ is a Borel graph. A particularly nice situation occurs when the action $\G \acts X$ is \emphd{free}, meaning that $\gamma \cdot x \neq x$ for all $x \in X$ and $\gamma \in \G \setminus \set{\mathbf{1}_\G}$. In this case, every connected component of $\mathsf{Sch}(X, F)$ is isomorphic to the Cayley graph $\mathsf{Cay}(\G, F)$. 
    
    A crucial example of a Borel action is the \emphd{\ep{Bernoulli} shift action} $\G \acts 2^\G$, given by the formula
    \[
        (\gamma \cdot x)(\delta) \,\defeq\, x(\delta \gamma) \quad \text{for all } x \colon \G \to 2 \text{ and } \gamma, \, \delta \in \G.
    \]
    Here we use the standard convention and identify $2$ with the two-element discrete space $\set{0,1}$. 
    The \emphd{free part} of $2^\G$, denoted by $\Free(2^\G)$, is the set of all points $x \in 2^\G$ whose stabilizer under the shift action is trivial. By definition, $\Free(2^\G)$ is the largest shift-invariant subspace of $2^\G$ on which the shift action is free. 
    When $\G$ is generated by a finite symmetric set $F \subseteq \G \setminus \set{\mathbf{1}_\G}$, we refer to the Schreier graph $\mathsf{Sch}(\Free(2^\G), F)$ as the \emphd{shift graph} of $\G$ and denote it by $\mathsf{Shift}(\G, F)$. 
    
    One of the most powerful approaches to studying finitely generated groups is to connect the algebraic properties of a group with the large-scale geometry of its Cayley graph. For example, a finitely generated group $\G$ is said to be \emphd{of polynomial growth} if its Cayley graph with respect to some \ep{equivalently, every} finite generating set is of polynomial growth. This ``geometric'' notion turns out to have an intimate relationship with the algebraic structure of $\G$; namely, by a celebrated theorem of Gromov \cite{gromov1981groups}, a finitely generated group is of polynomial growth if and only if it is virtually nilpotent.
    
    By Theorem~\ref{theo:KL}, if $\G$ is a group of polynomial growth, then its Cayley graph \ep{with respect to any finite generating set} admits an injective homomorphism to $\gr_{n,\infty}$ for some $n \in \N$, i.e., to the Cayley graph of the group $\Z^n$. We extend this result to Schreier graphs of such groups, and, even more generally, to arbitrary Borel graphs of polynomial growth. Namely, we show that if $G$ is a Borel graph with $\er(G) = \rho < \infty$, then $G$ admits an injective Borel homomorphism to the shift graph of $\Z^n$ corresponding to the generating set $\set{z \in \Z^n \,:\, \|z\|_\infty = 1}$, where $n = O(\rho \log \rho)$. Furthermore, we prove analogous extensions of all our other results \ep{i.e., the ones concerning asymptotically injective contractions and coarse embeddings}. As we explain in \S\ref{subsec:hyperfinite}, these results have a number of consequences in the study of Borel graphs; for instance, they imply that Borel graphs of polynomial growth are hyperfinite (see Definition \ref{defn:hyperfinite} and Corollary~\ref{corl:hyperfinite}).
    
    The principal difficulty in extending combinatorial results to the setting of Borel graphs arises from the need to ensure that all sets, functions, etc.\ that one considers are 
    \emph{Borel}.  
    This precludes working with each connected component of $G$ individually and requires treating them all in some sense ``uniformly.'' We outline some of the measures we take to deal with this challenge in \S\ref{subsec:roadmap}. 
    

    \subsection{Mapping graphs of polynomial growth to grid graphs}\label{subsec:results}
    
    After this preliminary discussion, we can now formally state our main results. For $n \in \N$, let
    \[
        \shgr_{n, \infty} \,\defeq\, \mathsf{Shift}(\Z^n, \set{z \in \Z^n \,:\, \|z\|_\infty = 1}).
    \]
    Note that every connected component of $\shgr_{n,\infty}$ is isomorphic to the $\infty$-grid graph $\gr_{n, \infty}$. For brevity, we write $\dist_\infty$ to denote the graph metric on $\shgr_{n,\infty}$ \ep{this notation is the same as our notation for the graph metric on $\gr_{n,\infty}$, but it will be clear from the context which of the two graphs we are referring to}. Explicitly, for two points $x$, $y \in \Free(2^{\Z^n})$, $\dist_\infty(x,y)$ is equal to the $\infty$-norm of the unique element $z \in \Z^n$ with $y = z \cdot x$ if such $z$ exists, and $\infty$ otherwise. 
    
    {
		\renewcommand{\arraystretch}{1.4}
		\begin{table}[t]
			\begin{tabular}{| l || c | c | c | c |}
				\hline
				& $b_\ell(r)$ & $b_u(r)$ & $n$ & injective? \\\hline
				Theorem~\ref{theo:coarse} & $\Omega_{G, \epsilon}(r^{1-\epsilon})$ & $r$ & $O_\epsilon(\ar)$ & no \\\hline
				Theorem~\ref{theo:injcoarse1} & $\Omega_{G, \epsilon}(r^{1-\epsilon})$ & $\max\set{r, O_{G, \epsilon}(1)}$ & $O_\epsilon(\ar)$ & yes \\\hline
				Theorem~\ref{theo:injcoarse2} & $\Omega_{\er, \epsilon}(r^{1-\epsilon})$ & $r$ & $O_\epsilon(\er \log \er)$ & yes \\\hline
			\end{tabular}
			\vspace*{10pt}
			\caption{The parameters of our embedding results. Here $n$ is the grid dimension and $b_\ell$, $b_u$ are the bounding functions for the coarse embedding \ep{as in Definition~\ref{defn:coarse}}.}\label{table:params}
		\end{table}	
	}
    Our central result yields a Borel coarse embedding from a Borel graph $G$ of polynomial growth to the graph $\shgr_{n,\infty}$ with $n$ linear in the asymptotic growth rate of $G$:
    
    \begin{theo}\label{theo:coarse}
        Let $G$ be a Borel graph with asymptotic growth rate $\ar(G) = \rho < \infty$. For every $\epsilon > 0$, there is a Borel map $f \colon V(G) \to \Free(2^{\Z^n})$, where $n = O_\epsilon(\rho)$, such that for all $u$, $v \in V(G)$,
        \[
            \dist_G(u, v) \,\geq\, \dist_\infty(f(u), f(v)) \,=\, \Omega_{G,\epsilon}(\dist_G(u, v)^{1-\epsilon}).
        \]
    \end{theo}
    
    In Theorem~\ref{theo:coarse_p}, we give a more precise (albeit more technical) statement that provides explicit bounds instead of the asymptotic notation used in Theorem~\ref{theo:coarse}. In particular, the bounds in Theorem~\ref{theo:coarse_p} are meaningful even if $G$ is a finite graph.
    
    Notice that if we apply Theorem~\ref{theo:coarse} to a \emph{connected} \ep{hence countable} graph $G$, then the image of $f$ would have to be contained in a single connected component of $\shgr_{n,\infty}$. Since every component of $\shgr_{n,\infty}$ is isomorphic to $\gr_{n,\infty}$, such $f$ provides a coarse embedding/asymptotically injective contraction from $G$ to $\gr_{n,\infty}$. For the benefit of the reader specifically interested in connected graphs, we record this observation here as a corollary:
    
    \begin{corl}\label{corl:coarse}
        Let $G$ be a connected graph with asymptotic growth rate $\ar(G) = \rho < \infty$. For every $\epsilon > 0$, there is a map $f \colon V(G) \to \Z^n$, where $n = O_\epsilon(\rho)$, such that for all $u$, $v \in V(G)$,
        \[
            \dist_G(u, v) \,\geq\, \dist_\infty(f(u), f(v)) \,=\, \Omega_{G,\epsilon}(\dist_G(u, v)^{1-\epsilon}).
        \]
    \end{corl}
    
    The function $f$ given by Theorem~\ref{theo:coarse} \ep{or Corollary~\ref{corl:coarse}} is both a coarse embedding \ep{this solves a problem of Papasoglu \cite[Question 4.1]{papasoglu2021polynomial}} and an asymptotically injective contraction. As mentioned earlier, the linear dependence of $n$ on $\rho$ can be seen as confirming part \ref{item:H2} of the Levin--Linial--London--Rabinovich conjecture ``in the asymptotic sense.''
    
    Next we consider the problem of constructing injective coarse embeddings. As advertised earlier, we prove two results in this direction, in the first of which we do not attempt to make $f$ a contraction, which allows us to keep $n$ linear as a function of $\ar(G)$: 
    
    \begin{theo}\label{theo:injcoarse1}
        Let $G$ be a Borel graph with $\ar(G) = \rho < \infty$. For every $\epsilon > 0$, there exists an injective Borel map $f \colon V(G) \to \Free(2^{\Z^n})$, where $n = O_\epsilon(\rho)$, such that for all $u$, $v \in V(G)$,
        \[
            \dist_\infty(f(u), f(v)) \,\leq\, \max \left\{\dist_G(u,v), \, O_{G,\epsilon}(1)\right\} \quad \text{and} \quad \dist_\infty(f(u), f(v)) \,=\, \Omega_{G,\epsilon}(\dist_G(u, v)^{1-\epsilon}).
        \]
    \end{theo}
    
    Again, by applying Theorem~\ref{theo:injcoarse1} to a connected graph $G$, we obtain the following corollary:
    
    \begin{corl}\label{corl:injcoarse1}
        Let $G$ be a connected graph with $\ar(G) = \rho < \infty$. For every $\epsilon > 0$, there exists an injective map $f \colon V(G) \to \Z^n$, where $n = O_\epsilon(\rho)$, such that for all $u$, $v \in V(G)$,
        \[
            \dist_\infty(f(u), f(v)) \,\leq\, \max \left\{\dist_G(u,v), \, O_{G,\epsilon}(1)\right\} \quad \text{and} \quad \dist_\infty(f(u), f(v)) \,=\, \Omega_{G,\epsilon}(\dist_G(u, v)^{1-\epsilon}).
        \]
    \end{corl}
    
    Finally, we prove a generalization of Theorem~\ref{theo:KL}\ref{item:KL1} that holds for Borel graphs and yields an injective homomorphism and a coarse embedding at the same time:
    
    \begin{theo}\label{theo:injcoarse2}
        Let $G$ be a Borel graph with $\er(G) = \rho < \infty$. For every $\epsilon > 0$, there exists an injective Borel map $f \colon V(G) \to \Free(2^{\Z^n})$, where $n = O_\epsilon(\rho \log \rho)$, such that for all $u$, $v \in V(G)$,
        \[
            \dist_G(u, v) \,\geq\, \dist_\infty(f(u), f(v)) \,=\, \Omega_{\rho,\epsilon}(\dist_G(u, v)^{1-\epsilon}).
        \]
    \end{theo}
    
    \begin{corl}\label{corl:injcoarse2}
        Let $G$ be a connected graph with $\er(G) = \rho < \infty$. For every $\epsilon > 0$, there exists an injective map $f \colon V(G) \to \Z^n$, where $n = O_\epsilon(\rho \log \rho)$, such that for all $u$, $v \in V(G)$,
        \[
            \dist_G(u, v) \,\geq\, \dist_\infty(f(u), f(v)) \,=\, \Omega_{\rho,\epsilon}(\dist_G(u, v)^{1-\epsilon}).
        \]
    \end{corl}
    
    The relationship between lower and upper bounds on $\dist_\infty(f(u), f(v))$ and the dimension of the grid in these results is summarized in Table \ref{table:params}. For numerically explicit versions of Theorems~\ref{theo:injcoarse1} and \ref{theo:injcoarse2}, see Theorem \ref{theo:injcoarse1_p}.

    \subsection{Applications: hyperfiniteness and toast}\label{subsec:hyperfinite}

    Our results stated in \S\ref{subsec:results} provide a general approach to questions about Borel graphs of polynomial growth: map the given graph $G$ to $\shgr_{n,\infty}$ for some $n \in \N$, solve the problem there by exploiting the grid structure, and then translate the solution back to $G$. In this subsection, we give some examples that illustrate the power of this method.

    A major source of motivation for our work is the theory of hyperfinite equivalence relations. As usual, we treat equivalence relations as sets of ordered pairs. In particular, an equivalence relation $E$ on a standard Borel space $X$ is said to be Borel if it is a Borel subset of $X \times X$. 
    
    \begin{defn}[Hyperfinite equivalence relations]\label{defn:hyperfinite}
        Let $E$ be a Borel equivalence relation on a standard Borel space $X$. We say that $E$ is \emphd{finite} if every $E$-class is finite, and \emphd{hyperfinite} if there exists an increasing sequence $E_0 \subseteq E_1 \subseteq E_2 \subseteq \ldots$ of finite Borel equivalence relations on $X$ whose union is $E$.
    \end{defn}
    
    Intuitively, Definition~\ref{defn:hyperfinite} says that a hyperfinite equivalence relation can be ``approximated'' by an increasing sequence of finite Borel equivalence relations. The systematic study of hyperfinite equivalence relations was initiated by Weiss \cite{Weiss} and Slaman and Steel \cite{SS}, with important foundational work done by Dougherty, Jackson, and Kechris \cite{DJK} and Jackson, Kechris, and Louveau \cite{jackson2002countable}, among others. For a survey of this subject, see \cite{KechrisCBER}.
    
    The Borel equivalence relations we shall consider come from graphs and group actions. Given a graph $G$, we let $\sim_G$ denote the \emphd{connectedness relation} of $G$, i.e., the equivalence relation on $V(G)$ whose classes are the connected components of $G$. We say that a Borel graph $G$ is \emphd{hyperfinite} if $\sim_G$ is a hyperfinite equivalence relation. Similarly, given a group action $\G \acts X$, we let $E(X, \G)$ be the corresponding \emphd{orbit equivalence relation}, i.e., the equivalence relation on $X$ whose classes are the orbits of the action, and we say that a Borel action $\G \acts X$ is \emphd{hyperfinite} if $E(X, \G)$ is a hyperfinite relation. Note that if $\G \acts X$ is a Borel action of a finitely generated group $\G$ on a standard Borel space, then its orbit equivalence relation coincides with the connectedness relation of the corresponding Schreier graph \ep{with respect to any finite generating set}.
    
    The following is one of the most notorious open questions in the area:
    
    \begin{ques}[{Weiss \cite{Weiss}}]\label{ques:Weiss}
        Is every Borel action of an amenable countable group on a standard Borel space hyperfinite?
    \end{ques}
    
    The assumption that $\G$ is amenable in Question~\ref{ques:Weiss} is necessary, since it is known that the shift action of a non-amenable countable group is not hyperfinite \cite[Corollary 1.8]{jackson2002countable}. Ornstein and Weiss \cite{ornstein1980ergodic} showed that the answer to Question~\ref{ques:Weiss} would be positive if we were allowed to discard a measure-$0$ subset of $V(G)$ with respect to some probability measure on $V(G)$. One may naively guess that handling the remaining measure-$0$ set should not pose much of a challenge, and yet Question~\ref{ques:Weiss} has remained wide open since the 1980s. 
    
    Until very recently, the affirmative answer to Question \ref{ques:Weiss} was known only in the case when $\G$ is of polynomial growth, due to a theorem of Jackson, Kechris, and Louveau \cite{jackson2002countable}:
    
    \begin{theo}[{Jackson--Kechris--Louveau \cite[Theorem 1.16]{jackson2002countable}}]\label{theo:JKL}
        If $\G$ is a group of polynomial growth, then every Borel action of $\G$ on a standard Borel space is hyperfinite.
    \end{theo}
    
    
    
    In a recent breakthrough \cite{conley2020borel}, Conley, Jackson, Marks, Seward, and Tucker-Drob  were able to give the affirmative answer to Question~\ref{ques:Weiss} for certain amenable groups of exponential growth, for example the lamplighter group $\Z_2 \wr \Z$, the Baumslag--Solitar group $BS(1,2)$, and all polycyclic groups. 
    In a different direction, it is natural to ask whether Theorem \ref{theo:JKL} holds for all Borel graphs of polynomial growth and not just for Schreier graphs: 
    
    \begin{ques}\label{ques:Marks}
        Is every Borel graph of polynomial growth hyperfinite?
    \end{ques}

    This is a well-known question in the area, whose importance was recently emphasized by Marks in his ICM address \cite[Problem 3.5]{MarksICM}. The methods used to establish Theorem \ref{theo:JKL} do not seem applicable to Question~\ref{ques:Marks}, because they crucially rely on having both upper and \emph{lower} bounds on the sizes of the balls in $G$, which are not available for general graphs of polynomial growth. To quote Marks \cite{MarksICM}:
    \begin{equation*}
    \parbox{0.87\textwidth}{``\emph{Finding techniques for resolving \eb{Question~\ref{ques:Marks}} where there is far less regular geometric structure would be one way of making progress towards resolving \hyperref[ques:Weiss]{Weiss's question}.}''}
    \end{equation*}

    We obtain a positive answer to Question~\ref{ques:Marks} as an immediate corollary of our main results:

    \begin{corl}\label{corl:hyperfinite}
        Every Borel graph of polynomial growth is hyperfinite.
    \end{corl}
    \begin{proof}
        Let $G$ be a Borel graph of polynomial growth. By using Theorem~\ref{theo:injcoarse2}, we may identify $G$ with its image under $f$ and assume that $G$ is a subgraph of $\shgr_{n,\infty}$ for some $n \in \N$. Since $\shgr_{n,\infty}$ is hyperfinite by Theorem~\ref{theo:JKL}, and since subgraphs of hyperfinite graphs are hyperfinite, we conclude that $G$ is hyperfinite.
    \end{proof}

    Using similar arguments, we can establish other properties of Borel graphs of polynomial growth by quoting known results about Schreier graphs of $\Z^n$-actions. One such property is the existence of a combinatorial structure called a \emph{toast}. It provides a particularly well-behaved witness to the hyperfiniteness of the associated connectedness relation and is very useful in combinatorial constructions. The definition of a toast stems from the work of Conley and Miller \cite{CM16} (the term ``toast'' was coined by Miller) and has been widely used in descriptive combinatorics \cite{CM16, bowen2021perfect, BPZ, trees, GLKS_Abelian, gao2022forcing, marks2017borel}. Related notions also appear in the theory of random processes \cite{BHT21,HSW17,Spi20}. Recently, Greb\'ik and Rozho\v{n} initiated a systematic study of toast-based constructions, under the name of \textsf{TOAST} algorithms \cite{GRgrids}.
    
    If $X$ is a standard Borel space, then the set $\fins{X}$ of all finite subsets of $X$ also carries a natural standard Borel structure; namely, a subset of $\fins{X}$ is Borel if and only if its preimage in $\bigsqcup_{k = 0}^\infty X^k$ under the map $(x_0, \ldots, x_{k-1}) \mapsto \set{x_0, \ldots, x_{k-1}}$ is Borel. When we call a family $\mathcal{F}$ of finite subsets of $X$ Borel, we simply mean that $\mathcal{F}$ is a Borel subset of $\fins{X}$.

    \begin{defn}[$r$-toast]\label{def: toast}
        Let $G$ be a Borel graph. 
        For $r \in \N$, a Borel family $\T \subseteq \fins{V(G)}$ of finite sets is an \emphdef{$r$-toast} if the following two conditions hold:
    \begin{enumerate}[label=\ep{\normalfont{}T\arabic*}]
        \item \label{item:T1} for every edge $uv \in E(G)$, there is some $K \in \T$ such that $u$, $v \in K$, and
        \item \label{item:T2} for distinct $K$, $L \in \T$, we have either $B_G(K, r) \cap L = \0$, $B_G(K,r) \subseteq L$, or $B_G(L,r) \subseteq K$. 
    \end{enumerate}
    \end{defn}

    We show that Borel graphs of polynomial growth admit an $r$-toast for any $r \in \N$:

    \begin{corl}\label{corl:toast}
        For every Borel graph $G$ of polynomial growth and every $r \in \N$, there exists an $r$-toast $\T \subseteq \fins{V(G)}$.
    \end{corl}
    \begin{proof}
        As in the proof of Corollary~\ref{corl:hyperfinite}, we may assume that $G$ is a subgraph of $\shgr_{n,\infty}$ for some $n < \infty$. By a result of Gao, Jackson, Krohne, and Seward (unpublished), there is an $r$-toast $\T^\ast \subseteq \fins{\Free(2^{\Z^n})}$ for $\shgr_{n, \infty}$ (see \cite[Appendix A]{marks2017borel} by Marks and Unger for a proof). It is now straightforward to verify that $\T \defeq \set{K \cap V(G) \,:\, K \in \T^\ast}$ is an $r$-toast for $G$. 
    \end{proof}

    Corollary~\ref{corl:toast} is a strengthening of Corollary~\ref{corl:hyperfinite}, since the existence of a $0$-toast implies hyperfiniteness \cite{FelixMeasure}. Another consequence of our main results is that Borel graphs of polynomial growth have finite Borel asymptotic dimension; this is also a strengthening of Corollary~\ref{corl:hyperfinite}. Since Borel asymptotic dimension plays a crucial role in the proofs of our main results, we discuss it in detail in the next subsection.
    
    \subsection{Further results: Borel asymptotic dimension}
    
    \subsubsection{Asymptotic dimension}
    
    One of the key tools that we use to establish our results in \S\ref{subsec:results} is the notion of asymptotic dimension, introduced by Gromov \cite[\S1.E]{gromov1993asymptotic}.
    
    \begin{defn}[Asymptotic dimension]\label{defn:ad}
        A family $\mathcal{U}$ of subsets of a metric space $(X,d)$ is:
        \begin{itemize}
            \item \emphd{uniformly bounded} if $\sup_{U \in \mathcal{U}} \diam U < \infty$;
            \item \emphd{$r$-disjoint} if $d(U, U') > r$ for all $U \neq U'$ in $\mathcal{U}$.
        \end{itemize}
        The \emphd{asymptotic dimension} of $X$, in symbols $\ad(X)$, is the minimum $n \in \N$ \ep{if it exists} such that for every $r > 0$, there are families $\mathcal{U}_0$, \ldots, $\mathcal{U}_n$ of subsets of $X$ with the following properties:
        \begin{itemize}
            \item each $\mathcal{U}_i$ is uniformly bounded and $r$-disjoint;
            \item $\mathcal{U} \defeq \bigcup_{i=0}^n \mathcal{U}_i$ is a cover of $X$, i.e., $\bigcup \mathcal{U} = X$.
        \end{itemize}
        If no such $n$ exists, then we set $\ad(X) \defeq \infty$.
    \end{defn}
    
    While Definition~\ref{defn:ad} is the original definition of asymptotic dimension, numerous equivalent definitions exist \cite[\S3.1]{bell2008asymptotic}. We will review some of them in \S\ref{sec:asdim}. 
    
    Gromov's motivation for introducing asymptotic dimension was to apply this concept to the study of finitely generated groups. A basic result is that the asymptotic dimension of a Cayley graph of a group $\G$ is independent of the choice of a finite generating set \cite[Corollary 51]{bell2008asymptotic}, so we are justified in referring to the asymptotic dimension $\ad(\G)$ of $\G$. See the survey \cite{bell2008asymptotic} by Bell and Dranishnikov for an overview of results about asymptotic dimension of groups.
    
    Unsurprisingly, we have $\ad(\Z^n) = n$ \cite[Example 2.2.6]{LSG}. Furthermore, all groups of polynomial growth have finite asymptotic dimension. Here is a proof sketch. By Gromov's theorem \cite{gromov1981groups}, groups of polynomial growth are virtually nilpotent. It is not hard to see that if $\G' \leq \G$ is a subgroup of $\G$ of finite index, then $\ad(\G') = \ad(\G)$ \cite[Corollary 55]{bell2008asymptotic}, so we have reduced the problem to the nilpotent case. The conclusion now follows since finitely generated nilpotent groups have finite asymptotic dimension by a theorem of Bell and Dranishnikov \cite[Corollary 9]{BD2006} \ep{their result holds, more generally, for polycyclic groups}.
    
    The argument sketched above relies on several deep group-theoretic results and cannot be generalized to graphs of polynomial growth that are not Cayley graphs. Hume \cite{Hume} observed that for vertex-transitive graphs, the problem can be reduced to the Cayley graph case using results of Sabidussi \cite{Sab} and Trofimov \cite{Tro}. Nevertheless, the question of whether all graphs of polynomial growth have finite asymptotic dimension remained open until quite recently. For example, this question was explicitly stated in \cite[1021]{vspakula2019relative} by \v{S}pakula and Tikuisis.
    
    Lately, asymptotic dimension attracted considerable attention of the graph theory community, leading to exciting new progress. In their breakthrough paper \cite{bonamy2021asymptotic}, Bonamy, Bousquet, Esperet, Groenland, Liu, Pirot, and Scott introduced powerful tools of structural graph theory to the study of asymptotic dimension and settled a number of open questions. In particular, they observed in \cite[Corollary 9.2]{bonamy2021asymptotic} that the Krauthgamer--Lee theorem \ep{Theorem \ref{theo:KL}} implies that a graph $G$ of exact growth rate $\rho < \infty$ satisfies $\ad(G) = O(\rho \log \rho)$ \ep{and hence the asymptotic dimension of $G$ is finite}. A more careful analysis of Krauthgamer and Lee's work yields the bound $\ad(G) = O(\rho)$; this is a consequence of \cite[Theorem 5.3]{krauthgamer2007intrinsic} and Corollary~\ref{corl:pdtoBad} below. The optimal bound on the asymptotic dimension of a graph in terms of its growth rate was obtained by Papasoglu \cite{papasoglu2021polynomial}: 
    
    \begin{theo}[{Papasoglu \cite{papasoglu2021polynomial}}]\label{theo:ad}
        Every graph $G$ satisfies $\ad(G) \leq \ar(G)$.
    \end{theo}
    
    Papasoglu's proof of Theorem~\ref{theo:ad} proceeds by induction on $\lfloor\ar(G)\rfloor$. As a byproduct of our investigation, we obtain an alternative, probabilistic proof of Theorem~\ref{theo:ad} based on some of the ideas involved in the Krauthgamer--Lee paper \cite{krauthgamer2007intrinsic}. Our proof is presented in \S\ref{sec:asdim}. Moreover, our proof approach also works in the setting of Borel graphs and yields a Borel version of Theorem~\ref{theo:ad}, to which we now turn our attention. 
    
    \subsubsection{Borel asymptotic dimension}
    
    Asymptotic dimension was introduced to descriptive combinatorics by Conley, Jackson, Marks, Seward, and Tucker-Drob \cite{conley2020borel}. 
    Recall that a graph is \emphd{locally finite} if all its vertices have finitely many neighbors. If $G$ is a locally finite Borel graph, then a subset of $V(G)$ of finite diameter in the graph metric must be finite. Thus, for such graphs we may formulate Definition~\ref{defn:ad} with the additional requirement that the families $\mathcal{U}_0$, \ldots, $\mathcal{U}_n$ be Borel subsets of $\fins{V(G)}$; this is exactly how Borel asymptotic dimension of $G$ is defined:
    
    \begin{defn}[Borel asymptotic dimension]\label{defn:Bad}
        Let $G$ be a locally finite Borel graph. The \emphd{Borel asymptotic dimension} of $G$, in symbols $\Bad(G)$, is the minimum $n \in \N$ \ep{if it exists} such that for every $r > 0$, there are Borel families $\mathcal{U}_0$, \ldots, $\mathcal{U}_n \subseteq \fins{V(G)}$ with the following properties:
        \begin{itemize}
            \item each $\mathcal{U}_i$ is uniformly bounded and $r$-disjoint (in the graph metric $\dist_G$);
            \item $\mathcal{U} \defeq \bigcup_{i=0}^n \mathcal{U}_i$ is a cover of $V(G)$.
        \end{itemize}
        If no such $n$ exists, then we set $\Bad(G) \defeq \infty$.
    \end{defn}
    
     The definition in \cite{conley2020borel} is slightly different from Definition~\ref{defn:Bad}, as it uses Borel equivalence relations in place of Borel families of finite sets; however, the two definitions are easily seen to be equivalent. 
     Also, the definition in \cite{conley2020borel} works for more general metric spaces, but we shall confine ourselves to the case of locally finite graphs for simplicity.
    
    Let $G$ be a locally finite Borel graph. Clearly, we have $\ad(G) \leq \Bad(G)$. A remarkable result of Conley \emph{et al.}\ is that this inequality can only be strict if $\Bad(G) = \infty$:
    
    \begin{theo}[{Conley--Jackson--Marks--Seward--Tucker-Drob \cite[Theorem 1.1]{conley2020borel}}]\label{theo:Bad_and_ad}
        Let $G$ be a locally finite Borel graph. If $\Bad(G) < \infty$, then $\Bad(G) = \ad(G)$.
    \end{theo}
    
    The work of Conley \emph{et al.}\ was motivated by applications of Borel asymptotic dimension to the theory of hyperfinite equivalence relations. 
    Specifically, they proved that finite Borel asymptotic dimension implies hyperfiniteness:
    
    \begin{theo}[{Conley--Jackson--Marks--Seward--Tucker-Drob \cite[Theorem 1.7]{conley2020borel}}]\label{theo:Bad_to_hyp}
        Let $G$ be a locally finite Borel graph. If $\Bad(G) < \infty$, then $G$ is hyperfinite.
    \end{theo}
    
    The utility of Theorem~\ref{theo:Bad_to_hyp} comes from the fact that for certain graphs, it is easier to bound Borel asymptotic dimension than to prove hyperfiniteness directly. Using this idea, Conley \emph{et al.}\ were able to provide the first examples of amenable groups of exponential growth for which the answer to Weiss's Question~\ref{ques:Weiss} is positive. 
    They also strengthened the \hyperref[theo:JKL]{Jackson--Kechris--Louveau Theorem} \ref{theo:JKL} as follows:
    
    \begin{theo}[{Conley--Jackson--Marks--Seward--Tucker-Drob \cite[Corollary 5.5]{conley2020borel}}]\label{theo:Bad_group}
        Let $\G$ be a group of polynomial growth. If $G$ is the Schreier graph of a Borel action of $\G$ on a standard Borel space, then $\Bad(G) < \infty$.
    \end{theo}
    
    
    
    
    We extend this result to all Borel graphs of polynomial growth and obtain the following 
    common generalization of Theorems~\ref{theo:ad} and \ref{theo:Bad_group}:
    
    \begin{theo}\label{theo:Bad}
        Every locally finite Borel graph $G$ satisfies $\Bad(G) \leq \ar(G)$. 
    \end{theo}
    
    In view of Theorems \ref{theo:ad} and \ref{theo:Bad_and_ad}, to establish Theorem~\ref{theo:Bad}, we only need to show that Borel graphs of polynomial growth have finite Borel asymptotic dimension. Since $\Bad(\shgr_{n,\infty}) < \infty$ by Theorem~\ref{theo:Bad_group}, and since Borel asymptotic dimension is monotone under taking subgraphs, Theorem~\ref{theo:Bad} is a consequence of Theorem~\ref{theo:injcoarse2} (analogous to Corollaries~\ref{corl:hyperfinite} and \ref{corl:toast}). However, we actually must verify Theorem~\ref{theo:Bad} without relying on Theorem~\ref{theo:injcoarse2}, since it (and some other, closely related facts) is used as a step in proving the results in \S\ref{subsec:results}. Our proof of Theorem~\ref{theo:Bad} follows a probabilistic approach inspired by the work of Krauthgamer and Lee \cite{krauthgamer2007intrinsic} and directly yields the bound $\Bad(G) \leq \ar(G)$ (thus giving an alternative proof of Theorem~\ref{theo:ad}). In particular, since finite Borel asymptotic dimension implies hyperfiniteness, our proof of Theorem~\ref{theo:Bad} presented in \S\ref{sec:asdim} also gives a streamlined proof of Corollary~\ref{corl:hyperfinite} that does not require establishing the (considerably more involved) results in \S\ref{subsec:results}. We anticipate that the techniques we develop to prove Theorem~\ref{theo:Bad} will have further applications to the theory of hyperfiniteness. 
    
    \subsection{A road map}\label{subsec:roadmap}
    
    Let us now outline the structure of the remainder of the paper. We should emphasize that, although all our arguments are presented in the Borel setting, a reader who is only interested in connected graphs may simply ignore every appearance of the word ``Borel'' in the remainder of the paper and mentally replace the words ``standard Borel space'' by ``set'' everywhere.
    
    \subsubsection*{\normalfont{\itshape{}The Lov\'asz Local Lemma}}
    
    The proofs of our main results have a significant probabilistic component. Specifically, they rely in a crucial way on the so-called \emphd{Lov\'asz Local Lemma}, or the \emphd{LLL} for short. The LLL was introduced by Erd\H{o}s and Lov\'asz in the 1970s \cite{EL} and has by now become an indispensable tool throughout  combinatorics. The LLL is typically used to prove the existence of colorings and other combinatorial objects satisfying given sets of ``local'' constraints. As with other existence results, it is a matter of interest to determine whether the LLL can be used to derive conclusions that are ``constructive'' in various senses. The first such ``constructive'' version of the LLL was the algorithmic LLL due to Beck \cite{Beck}. Beck's result requires somewhat stronger numerical assumptions than the ordinary LLL. This discrepancy has been eventually eliminated in the breakthrough work of Moser and Tardos \cite{MT}. The Moser--Tardos method was later adapted to derive ``constructive'' analogs of the LLL in a variety of different contexts. Most pertinently for our present work, Cs\'{o}ka, Grabowski, M\'{a}th\'{e}, Pikhurko, and Tyros used a variant of the Moser--Tardos technique to establish a Borel version of the LLL for graphs of subexponential growth \cite{CGMPT}. Since we are working with graphs of polynomial growth in this paper, we are able to make full use of the Cs\'{o}ka \emph{et al.}\ Borel LLL in our arguments. We introduce the necessary definitions and state the Borel LLL in \S\ref{sec:LLL}.
    
    \subsubsection*{\normalfont{\itshape{}Padded decompositions and asymptotic dimension}}
    
    We establish Theorem~\ref{theo:Bad} \ep{i.e., a bound on the Borel asymptotic dimension of graphs of polynomial growth} and some related results in \S\ref{sec:asdim}. Our arguments in this section are directly inspired by the proof of \cite[Theorem 5.3]{krauthgamer2007intrinsic} due to Krauthgamer and Lee, which is in turn based on the work of Linial and Saks \cite{LinialSaks} and Bartal \cite{Bartal}. Instead of Definition~\ref{defn:Bad}, we use an equivalent definition of \ep{Borel} asymptotic dimension that involves so-called \emph{padded decompositions} of graphs. Roughly speaking, a padded decomposition of a graph $G$ is a family of partitions of $V(G)$ into finite clusters such that every vertex is far from the boundary of its cluster in at least one of the partitions; the details are given in \S\ref{sec:asdim}. In order to prove Theorem~\ref{theo:Bad}, we construct a family of partitions via a randomized procedure and then invoke the Borel LLL to argue that this procedure can successfully yield a padded decomposition with the desired properties. We then apply a similar procedure to obtain padded decompositions with a different range of parameters, which play an important role in our proof of Theorem~\ref{theo:coarse}. While the construction we use is similar to the one in \cite{krauthgamer2007intrinsic} \ep{except that it is preformed in the Borel setting}, our analysis is significantly more precise and more intricate, because we require much better control over the resulting numerical parameters.
    
    \subsubsection*{\normalfont{\itshape{}Cocycles}}
    
    Our embedding results in \S\ref{subsec:results} assert the existence of functions of the form $f \colon V(G) \to \mathsf{Free}(2^{\Z^n})$ with some properties. Instead of building such $f$ directly, it turns out to be more convenient to work with certain auxiliary structures, called \emph{cocycles}. Recall that $\sim_G$ denotes the connectedness relation of $G$, which we think of as a set of ordered pairs of vertices, i.e., a subset of $V(G) \times V(G)$. Assume that $f \colon V(G) \to \mathsf{Free}(2^{\Z^n})$ is a mapping such that $\dist_\infty(f(u), f(v)) < \infty$ whenever $\dist_G(u,v) < \infty$ (for instance, if $f$ is a coarse embedding of $G$ into $\mathsf{ShiftGrid}_{n, \infty}$, it must have this property). We can then define a function $\updelta f \colon {\sim_G} \to \Z^n$, called the cocycle associated to $f$, as follows. Let $u$ and $v$ be vertices such that $u \sim_G v$, i.e., $(u,v) \in {\sim_G}$. Then $\dist_G(u,v) < \infty$ and so $\dist_\infty(f(u), f(v)) < \infty$ as well, which means that $f(u)$ and $f(v)$ belong to the same orbit of the shift action $\Z^n \acts \Free(2^{\Z^n})$. Therefore, we can let $\updelta f(u,v) \in \Z^n$ be the unique group element such that \[(\updelta f (u,v)) \cdot f(u) \,=\, f(v).\] In \S\ref{sec:cocycle}, we develop tools that allow us to reverse this construction. In other words, given a function $\updelta \colon {\sim_G} \to \Z^n$ with certain properties, we are able to find a mapping $f$ such that $\updelta = \updelta f$. The main advantage of working with cocycles instead of dealing with maps from $V(G)$ to $\Free(2^{\Z^n})$ directly is that the sum of two cocycles is again a cocycle. This linearity enables an inductive approach wherein a cocycle with the desired properties is built in stages as a sum of several pieces.
    
    \subsubsection*{\normalfont{\itshape{}Coarse embeddings}}
    
    In \S\ref{sec:coarse} we prove our main Theorem~\ref{theo:coarse}, asserting the existence of a Borel coarse embedding from a Borel graph of polynomial growth to the shift graph of $\Z^n$. The overall proof strategy we employ is analogous to the proof of the Krauthgamer--Lee Theorem~\ref{theo:KL} \cite{krauthgamer2007intrinsic} and is inspired by the earlier work of Rao \cite{rao1999}, although we streamline and sharpen many parts of the argument. The proof utilizes the padded decompositions built in \S\ref{sec:asdim} and proceeds in countably many stages, where on each stage we ensure that the desired relation $\dist_\infty(f(u), f(v)) = \Omega_{G,\epsilon}(\dist_G(u, v)^{1-\epsilon})$ holds for a certain small range of values of $\dist_G(u,v)$. We have to take great care to ensure that together these ranges cover all sufficiently large positive real numbers. Each step of the construction is probabilistic and employs the Borel LLL. In particular, we use the polynomial growth assumption on two independent occasions in the proof: once to construct the requisite padded decompositions and the second time to create the coarse embedding itself. The cocycle machinery developed in \S\ref{sec:cocycle} allows us to ``glue'' the outcomes of the individual steps together in the end.
    
    \subsubsection*{\normalfont{\itshape{}Making the maps injective}}
    
    Finally, in \S\ref{sec:injective}, we use Theorem~\ref{theo:coarse} \ep{or rather its numerically explicit version, Theorem~\ref{theo:coarse_p}} to derive Theorems~\ref{theo:injcoarse1} and \ref{theo:injcoarse2}. Here we again employ the cocycle machinery of \S\ref{sec:cocycle}: it turns out that an asymptotically injective function can be made truly injective fairly easily by a suitable modification of its associated cocycle.
    
    \subsubsection*{Acknowledgments}
    
    We are grateful to Andrew Marks, Alex Kastner, and Forte Shinko for helpful discussions during the preparation of this manuscript and to Xinran Tao for the picture of a dumpling. We also thank the anonymous referee for reading the manuscript carefully and providing helpful comments.

    \section{The Lov\'asz Local Lemma and its Borel version}\label{sec:LLL}
    \subsection{Constraint satisfaction problems and the LLL}
    
    In this paper we employ the Lov\'asz Local Lemma in the setting of constraint satisfaction problems. The following notation and terminology are borrowed from \cite{BerDist}.
    
    \begin{defn}
    Fix a set $X$ and a finite set $C$. A $C$-\emphdef{coloring} of a set $S$ is a function $f \colon S \rightarrow C$.
        \begin{itemize}[wide]
     \item (Constraints) Given a finite subset $D \subseteq X$, an \emphdef{$(X, C)$-constraint} (or simply a \emphdef{constraint} if $X$ and $C$ are clear from the context) with \emphdef{domain} $D$ is a set $A \subseteq C^D$ of $C$-colorings of $D$. We write $\dom(A) \defeq D$. A $C$-coloring $f \colon X \rightarrow C$ of $X$ \emphdef{violates} a constraint $A$ with domain $D$ if the restriction of $f$ to $D$ is in $A$, and \emphdef{satisfies} $A$ otherwise.
     \item (Constraint satisfaction problems) A \emphdef{constraint satisfaction problem} (or a \emphdef{CSP} for short) $\B$ on $X$ with range $C$, in symbols $\B \colon X \rightarrow^{?} C$, is a set of $(X, C)$-constraints. A \emphdef{solution} to a CSP $\B \colon X \rightarrow^{?} C$ is a $C$-coloring $f \colon X\rightarrow C$ that satisfies every constraint $A \in \B$.
     \end{itemize}
     \end{defn}

    In other words, each constraint $A \in \B$ in a CSP $\B \colon X \to^? C$ is viewed as a set of finite ``forbidden patterns'' that may not appear in a solution $f \colon X \to C$. The LLL provides a probabilistic condition that is sufficient to ensure that a given CSP has a solution. To state this condition, we need some further notation. Let $\B \colon X \rightarrow^{?} C$ be a CSP. For each constraint $A \in \B$, the \emphdef{probability} $\P[A]$ of $A$ is defined as the probability that $A$ is violated by a uniformly random $C$-coloring $f \colon X \to C$, that is,
 \[\P[A] \,\defeq\, \frac{|A|}{|C|^{|\dom(A)|}}.\]
 The \emphdef{neighborhood} of $A$ in $\B$ is the set 
 \[N(A) \,\defeq\, \set{A' \in \B \,:\, A' \neq A  \text{ and } \dom(A') \cap \dom(A) \neq \0}.\]
 Let $\p(\B) \defeq \sup_{A \in \B} \P[A]$ and $\d(\B) \defeq \sup_{A \in \B} |N(A)|$.
 

\begin{theo}[{Lov\'asz Local Lemma \cite{EL, SpencerLLL}}]
\label{thm: LLL}
If $\B$ is a CSP such that $e \cdot \p(\B) \cdot (\d(\B) + 1) < 1$,
then $\B$ has a solution. \ep{Here $e = 2.71\ldots$ is the base of the natural logarithm.}
\end{theo}

    The LLL is often stated for finite $\B$. However, a straightforward compactness argument shows that Theorem \ref{thm: LLL} holds for infinite $\B$ as well (see, e.g., the proof of \cite[Theorem 5.2.2]{AS}).

\subsection{Borel colorings and the Borel LLL}
 
In descriptive combinatorics, we are often interested in Borel solutions to CSPs defined on standard Borel spaces. Let $X$ be a standard Borel space and let $C$ be a finite set, which we view as a discrete space. Then every $(X, C)$-constraint $A$ can be viewed as a finite subset of $\fins{X \times C}$ by identifying each function $\phi \in A$ with its graph $\set{(x, \phi(x)) \,:\, x \in \dom(A)}$. Hence we may speak of \emphdef{Borel CSPs} $\B \colon X \rightarrow^{?} C$, i.e., Borel sets $\B \subseteq \fins{\fins{X \times C}}$ of $(X, C)$-constraints.

    In general, a Borel CSP $\B$ satisfying the conditions of Theorem~\ref{thm: LLL} may not have a Borel solution \cite[Theorem~1.6]{CJMST-D}. Nevertheless, Cs\'{o}ka, Grabowski, M\'{a}th\'{e}, Pikhurko, and Tyros were able to prove a Borel version of the LLL under an additional subexponential growth assumption on $\B$ \cite{CGMPT}. To state their result, we need some terminology. Given a CSP $\B \colon X \rightarrow^? C$, we define a graph $G_{\B}$ with vertex set $X$ by making two vertices $x \neq y$ adjacent if and only if there is a constraint $A \in \B$ such that $\set{x, y} \subseteq \dom(A)$. 
    A graph $G$ is of \emphdef{subexponential growth} if $\gamma_G(r) = e^{o(r)}$, i.e., if for every $\epsilon > 0$ there is $r > 0$ such that for all $R \geq r$ and all $v \in V(G)$ we have $|B_G(v, R)| \leq (1+\epsilon)^R$. Note that if $G$ is of polynomial growth, then it is of subexponential growth as well. 

\begin{theo}[{Borel Lov\'asz Local Lemma, Cs\'{o}ka--Grabowski--M\'{a}th\'{e}--Pikhurko--Tyros \cite[Theorem 4.5]{CGMPT}}]
\label{thm: BLL}
Let $\B \colon X \rightarrow^{?} C$ be a Borel CSP on a standard Borel space $X$. If $G_{\B}$ is of subexponential growth and 
$e \cdot \p(\B) \cdot (\d(\B) + 1) < 1$,
then $\B$ has a Borel solution $f \colon X \to C$. 
\end{theo}
    We remark that \cite[Theorem 4.5]{CGMPT} uses slightly different terminology from ours. However, Theorem~\ref{thm: BLL} is easily seen to be a consequence of \cite[Theorem 4.5]{CGMPT}; see \cite[Remark 1.3]{CGMPT}. 

\begin{remk}\label{remk:LLL}
    In the proof of Theorem~\ref{theo:Bad}, 
    we will need to use the Borel LLL in a somewhat more general setting than the one we have just described. 
    Namely, we shall consider the case when the set $C$ is equipped with a certain non-uniform probability distribution. Thankfully, we can reduce this more general situation to the one described above as follows. 

    Let us denote the given probability distribution by $\mu$, so $\mu(c)$ is the probability assigned to each $c \in C$. In this case it makes sense to define the probability of a constraint $A$ by
    \begin{equation}\label{eq:nonuni}
        \P[A] \,\defeq\, \sum_{\phi \in A} \prod_{x \in \dom(A)} \mu(\phi(x)).
    \end{equation}
    In other words, $\P[A]$ is the probability that $A$ is violated by a random coloring $f \colon X \to C$ where each color $f(x)$ is drawn independently from the distribution $\mu$.

    Suppose that $\B \colon X \to^? C$ is a Borel CSP on a standard Borel space $X$ satisfying the assumptions of Theorem~\ref{thm: BLL}, where the probabilities of constraints are defined via \eqref{eq:nonuni}. The graph $G_\B$ is of subexponential growth, so in particular its maximum degree is finite. It follows that there is some $s \in \N$ such that $|\dom(A)| \leq s$ for all $A \in \B$. Therefore, if we modify $\mu$ by changing each value $\mu(c)$ at most by a factor of $1 + \epsilon$ for some $\epsilon > 0$, then for all $A \in \B$, the value $\P[A]$ will change at most by a factor of $(1+\epsilon)^s$. Thus, since the inequality in Theorem~\ref{thm: BLL} is strict, we may slightly perturb $\mu$ and assume that $\mu(c)$ is rational for all $c \in C$.
    
    Now let us write $\mu(c) = n_c/d_c$, where $n_c$, $d_c \in \N$. Let $d \defeq \prod_{c \in C} d_c$ and $D \defeq \set{1, 2, \ldots, d}$. Fix a function $h \colon D \to C$ such that $|h^{-1}(c)| = \mu(c) d$ for all $c \in C$. Then, instead of picking $f(x)$ from the distribution $\mu$, we may choose $q(x) \in D$ uniformly at random and set $f(x) \defeq h(q(x))$. For every constraint $A \in \B$, there is a corresponding $(X, D)$-constraint $A'$ given by 
\[A' \,\defeq\, \{\phi \colon \dom(A) \to D \,:\, \textrm{the map } x \mapsto h(\phi(x)) \textrm{ is in  }\B\},\]
    and, clearly, $\P[A] = \P[A']$. In this way, we replace $\B$ by an ``equivalent'' CSP $\B' \defeq \{A' \,:\, A \in \B\}$ with range $D$. Note that $G_\B = G_{\B'}$ and hence $\B'$ also satisfies the assumptions of Theorem~\ref{thm: BLL}. Furthermore, if $q \colon X \to D$ is a Borel solution to $\B'$, then $x \mapsto h(q(x))$ is a Borel solution to $\B$. Therefore, we may conclude from Theorem~\ref{thm: BLL} that the original CSP $\B$ admits a Borel solution.
\end{remk}

\section{Padded decompositions and asymptotic dimension}\label{sec:asdim}

\subsection{Some background facts}

    The aim of this section is to prove Theorem~\ref{theo:Bad}, i.e., to show that a locally finite Borel graph $G$ satisfies $\Bad(G) \leq \ar(G)$. We also establish some related results that will be used in the proof of Theorem~\ref{theo:coarse}.  


    
    Before we proceed, let us briefly review some basic facts from descriptive set theory. As a reminder, our references for descriptive set theory are \cite{KechrisDST} by Kechris and \cite{AnushDST} by Tserunyan. 
    Throughout the paper, we shall use the following deep result without mention:
    
    \begin{theo}[{Luzin--Novikov \cite[Theorem 18.10]{KechrisDST}}]\label{theo:LN}
        Let $X$, $Y$ be standard Borel spaces and let $R \subseteq X \times Y$ be a Borel subset. Suppose that for all $x \in X$, the set $\set{y \in Y \,:\, (x,y) \in R}$ is countable. Then there exists a sequence $(f_n)_{n \in \N}$ of Borel partial functions $f_n \colon X \dashrightarrow Y$ defined on Borel subsets of $X$ such that $(x,y) \in R$ if and only if $y = f_n(x)$ for some $n \in \N$.
    \end{theo}
    
    This theorem is important because it allows defining Borel sets by quantifying over countable sets. For example, suppose that $G$ is a locally countable Borel graph (meaning that every vertex of $G$ has countably many neighbors). Then the connectedness relation $\sim_G$, viewed as a set of ordered pairs, is a Borel subset of $V(G) \times V(G)$. Indeed, we have
    \begin{align}
        u \sim_G v \quad \Longleftrightarrow \quad \exists d \in \N \, \exists (u_0, &\ldots, u_d) \in V(G)^{d+1}, \nonumber\\
        &u = u_0, \ u_0 u_1 \in E(G), \ \ldots,\ u_{d-1} u_d \in E(G), \ u_d = v.\label{eq:example}
    \end{align}
    Note that, since $G$ is locally countable, for any given pair of vertices $u$, $v$, there are only countably many tuples $(u_0, \ldots, u_d)$ satisfying the condition in \eqref{eq:example}. A standard argument using the \hyperref[theo:LN]{Luzin--Novikov theorem} then shows that the set $\sim_G$ is Borel. Explicitly, for $d \in \N$, define
    \[
        R_d \,\defeq\, \set{(u,v, u_0, \ldots, u_d) \in V(G)^{d+3} \,:\, u = u_0, \ u_0 u_1 \in E(G), \ \ldots,\ u_{d-1} u_d \in E(G), \ u_d = v}.
    \]
    The set $R_d$ is clearly Borel since $G$ is Borel. 
    Hence, by Theorem~\ref{theo:LN}, we have a sequence of Borel partial functions $f_{d,n} \colon V(G)^2 \dashrightarrow V(G)^{d+1}$ such that
    \[
        (u,v, u_0, \ldots, u_d) \in R_d \quad \Longleftrightarrow \quad \exists n \in \N, \ f_{d,n}(u,v) = (u_0, \ldots, u_d).
    \]
    Therefore, we can write
    \[
        u \sim_G v \quad \Longleftrightarrow \quad \exists d \in \N \, \exists n \in \N, \ (u,v) \in \dom(f_{d,n}),
    \]
    which means that $\sim_G$ is a countable union of Borel sets, so it is itself Borel. A similar argument shows that, for example, the function $\dist_G \colon V(G) \times V(G) \to \N \cup \set{\infty}$ is Borel. See \cite[\S5.2]{pikhurko2020borel} for a further discussion. Since the arguments analogous to the above are standard, they shall be omitted in the rest of the article.
    
    An important corollary of the \hyperref[theo:LN]{Luzin--Novikov theorem} is that countable-to-one images of Borel sets under Borel functions are Borel:
    
    \begin{corl}[{\cite[Corollary 13.7]{AnushDST}}]\label{corl:ctbl-to-one}
        Let $X$, $Y$ be standard Borel spaces and let $f \colon X \to Y$ be a Borel function. If $f$ is countable-to-one, then the image $f(X)$ of $f$ is a Borel subset of $Y$.
    \end{corl}
    
    Let $G$ be a graph. A \emphd{proper coloring} of $G$ is a function $c$ defined on $V(G)$ such that $c(u) \neq c(v)$ for all $uv \in E(G)$. The \emphd{Borel chromatic number} of a Borel graph $G$ is the smallest cardinality of a standard Borel space $Y$ such that there is a Borel proper coloring $c \colon V(G) \to Y$. This notion was introduced in the seminal paper \cite{kechris1999borel} by Kechris, Solecki, and Todorcevic. Among other results, they proved the following:
    
    \begin{theo}[{\cite[Proposition 4.6]{kechris1999borel}}]\label{theo:Delta+1}
        If $G$ is a Borel graph of finite maximum degree $\Delta$, then the Borel chromatic number of $G$ is at most $\Delta + 1$.
    \end{theo}
    
    A set $I \subseteq V(G)$ is \emphd{independent} in $G$ if $uv \not \in E(G)$ for all $u$, $v \in I$. Note that if $c \colon V(G) \to Y$ is a Borel proper coloring of $G$, where $|Y| = k$, then $V(G) = \bigcup_{y \in Y} c^{-1}(y)$ is a partition of $V(G)$ into $k$ Borel independent sets.

\subsection{Padded decompositions}

    As mentioned in \S\ref{subsec:roadmap}, instead of Definition~\ref{defn:Bad}, it will be more convenient for our purposes to use an equivalent way of describing \ep{Borel} asymptotic dimension in terms of auxiliary structures called \emph{padded decompositions}. We start with a few definitions.


 
 \begin{defn}[Basic terminology]
    Let $G$ be a graph.
    The \emphdef{\ep{inner} boundary} of a set $S \subseteq V(G)$ is \[\partial S \defeq \{u \in S \,:\, \exists uv \in E(G) \text{ with } v \notin S\}.\]
    The \emphdef{diameter} of $S$ is defined to be $\diam S \defeq \sup_{u, v  \in S} \dist_G(u, v)$. A family $\mathcal{U}$ of subsets of $V(G)$ is called \emphd{$D$-bounded} for some $D > 0$ if $\diam S \leq D$ for all $S \in \mathcal{U}$.
\end{defn}

\begin{defn}[Partitions and clusters]
    A \emphd{partition} of a set $X$ is a collection $\mP$ of pairwise disjoint subsets of $X$, called \emphd{clusters}, whose union is $X$. Given a set $S \subseteq X$ and a partition $\mP$ of $X$, 
    the corresponding \emphdef{induced partition} of $S$ is $S/\mathcal{P} \defeq \{C \cap S \,:\, C \in \mP\}$.
 \end{defn}
   We adopt the term ``clusters'' to be consistent with the terminology in Krauthgamer and Lee's paper \cite{krauthgamer2007intrinsic}. We remark that if $G$ is a graph and $\mP$ is a partition of $V(G)$, the clusters of $\mP$ 
   need not be connected.

\begin{defn}[$(r,D)$-covers]
\label{def:cover}
    For any graph $G$, a tuple $(\mU_{1}, \mU_{2}, \ldots, \mU_{m})$ of $m$ families of subsets of $V(G)$ is called 
     an \emphdef{$(r, D)$-cover} of $G$ with $m$ \emphdef{layers} if the following properties are satisfied. 
    \begin{enumerate}[label=\ep{\normalfont\arabic*}]
        \item\label{item:cover:rD} Each $\mU_i$ is $r$-disjoint and $D$-bounded (cf.\ Definition~\ref{defn:ad}).  
        \item\label{item:cover:cover} $\bigcup_{i=1 }^m \mU_i$ is a cover of $V(G)$.
    \end{enumerate}
\end{defn}
    As explained in the paragraph immediately preceding Definition~\ref{def: toast}, we call a family $\mU$ of finite subsets of a standard Borel space $X$ Borel if it is a Borel subset of the space $\fins{X}$. When we say that $\mU$ is a Borel family of subsets of $X$, it is tacitly implied that the sets in $\mU$ are finite. We say that a tuple $(\mU_1, \ldots, \mU_m)$ of families of subsets of $X$ is Borel if every $\mU_i$ is Borel. Note that if $G$ is a locally finite Borel graph and $\mU$ is a uniformly bounded Borel family of subsets of $V(G)$, then every vertex $v \in V(G)$ belongs to finitely many sets in $\mU$. Therefore, the map
    \[
        \set{(U, v) \,:\, U \in \mU, \, v \in U} \to V(G) \quad \colon \quad (U, v) \mapsto v
    \]
    is finite-to-one. By Corollary~\ref{corl:ctbl-to-one}, it follows that the image of this map, i.e., the set $\bigcup \mU$, is a Borel subset of $V(G)$. The following lemma connects the notion of a Borel partition of a locally finite graph to the concept of a Borel equivalence relation:
    
    \begin{Lemma}
        Let $G$ be a locally finite Borel graph. Suppose that $\mP$ is a partition of $V(G)$ such that each cluster $C \in \mP$ has finite diameter in $G$. Then the following statements are equivalent:
        \begin{itemize}
            \item $\mP$ is Borel;
            \item the equivalence relation $E_{\mP}$ on $V(G)$ whose classes are the $\mP$-clusters is Borel.
        \end{itemize}
    \end{Lemma}
    \begin{proof}
        Suppose $\mP$ is Borel. Then we  write
        \[
            u \, E_{\mP} \, v \quad \Longleftrightarrow \quad \exists C \in \mP, \ u \in C \text{ and } v \in C.
        \]
        Since there is only one cluster $C$ containing any given vertex $u$, we conclude the set $E_{\mP}$ is Borel by the \hyperref[theo:LN]{Luzin--Novikov theorem} (Theorem~\ref{theo:LN}). Conversely, suppose that $E_{\mP}$ is Borel. Then for any tuple $(u_0, \ldots, u_{k-1}) \in V(G)^k$, we have
        \begin{align*}
            \set{u_0, \ldots, u_{k-1}} \in \mP& \quad \Longleftrightarrow \quad \left(\forall i, j \in \set{0, \ldots, k-1}, \ u_i \, E_{\mP} \, u_j\right) \\
            &\text{and } \forall u \in V(G), \ u \, E_\mP \, u_0 \, \Longrightarrow \, (u = u_i \text{ for some $i \in \set{0, \ldots, k-1}$}).
        \end{align*}
        Since for any $u_0 \in V(G)$, there are only finitely many vertices $u \in V(G)$ such that $u \,E_\mP\, u_0$, it follows that $\mP$ is Borel by the \hyperref[theo:LN]{Luzin--Novikov theorem}.
    \end{proof}
    
    It follows from Definitions~\ref{defn:ad}, \ref{defn:Bad}, and \ref{def:cover} that the \ep{Borel} asymptotic dimension of a locally finite \ep{Borel} graph $G$ is the smallest $n \in \N$ \ep{if it exists} such that for every $r > 0$, there is $D > 0$ such that $G$ admits a \ep{Borel} $(r, D)$-cover with $n + 1$ layers.

\begin{defn}[Padded decompositions]
\label{def:padded_decomposition}
    Let $G$ be a graph and fix parameters $\alpha > 1$ and $r > 0$. 
    A tuple $(\mP_{1}, \mP_{2}, \ldots, \mP_{m})$ of $m$ partitions of $V(G)$ is called an \emphdef{$(r, \alpha)$-padded decomposition} of $G$ with $m$ \emphdef{layers} if the following properties are satisfied:
    \begin{enumerate}[label=\ep{\normalfont\arabic*}]
        \item\label{item:pd:diam} Each $\mathcal{P}_i$ is $r^\alpha$-bounded. 
        \item\label{item:pd:padded} For every $v \in V(G)$ there exists some partition $\mP_i$ in which there is $C \in \mP_i$ with $B_G(v, r) \subseteq C$. 
    \end{enumerate}
\end{defn}

    The relationship between Definitions~\ref{def:cover} and \ref{def:padded_decomposition} is given by the following lemma: 
    
    
\begin{Lemma}
\label{lem: TFAE}
    Let $G$ be a locally finite Borel graph. Then:
    \begin{enumerate}[label=\ep{\normalfont\arabic*}]
        \item\label{item:from_cover_to_pd} If there exists a Borel $(2r, D)$-cover of $G$ with $m$ layers, then there exists a Borel $(r, \alpha)$-padded decomposition of $G$ with $m$ layers for any $\alpha$ such that $r^{\alpha} \ge D+2r$.  
        \item\label{item:from_pd_to_cover} If there exists a Borel $(r/2 + 1, \alpha)$-padded decomposition of $G$ with $m$ layers, then there exists a Borel $(r, D)$-cover of $G$ with $m$ layers for any $D \ge (r/2 + 1)^{\alpha}$.
    \end{enumerate}
\end{Lemma}

\begin{proof}
    \ref{item:from_cover_to_pd} Let $(\mU_1, \mU_2, \ldots, \mU_m)$ be a Borel $(2r, D)$-cover of $G$ with $m$ layers for some $D > 0$. For each $i$, we modify $\mU_i$ to construct a partition $\mP_i$. Namely, we ``expand'' every $C \in \mU_i$ to $C' \defeq B_G(C,r)$. 
    Next we let $V_i \defeq \bigcup \{C' \,:\, C \in \mU_i\}$ and define
    \[
        \mP_i \,\defeq\, \{C' \,:\, C \in \mU_i\} \cup \{\{x\} \,:\, x \in V(G) \setminus V_i\}.
    \]
    Then $\mP_i$ is a Borel partition of $V(G)$. Pick $\alpha$ such that $r^{\alpha} \ge D+2r$. We claim that $(\mP_{1}, \mP_{2}, \ldots, \mP_{m})$ is a Borel $(r, \alpha)$-padded decomposition of $G$. 
Condition \ref{item:pd:diam} of Definition~\ref{def:padded_decomposition} is satisfied by the choice of $\alpha$.   To verify condition \ref{item:pd:padded}, note that 
for each $v \in V(G)$, there is some $C \in \mU_i$ such that $v \in C$, 
which implies that $B_G(v, r) \subseteq C'$, as desired. 


    \ref{item:from_pd_to_cover} Let $(\mP_{1}, \mP_{2}, \ldots, \mP_{m})$ be a $(r/2 + 1, \alpha)$-padded decomposition of $G$ for some $\alpha > 1$. 
For each $i$, we modify $\mP_i$ to construct a family of sets $\mU_i$. Namely, for every cluster $C \in \mP_i$, we ``shrink'' it to \[C' \,\defeq\, C\,\setminus \bigcup_{x \in \partial C} B_G(x, r/2).\] Define $\mU_i \defeq \{C' \,:\, C \in \mP_i\}$ and pick $D$ so that $D \ge (r/2 + 1)^{\alpha}$. We claim that $(\mU_1, \mU_2, \ldots, \mU_m)$ is a Borel $(r, D)$-cover of $G$ with $m$ layers. 
Condition \ref{item:cover:rD} of Definition~\ref{def:cover} is satisfied because the clusters in $\mP_i$ are disjoint and by the choice of $D$. 
To verify \ref{item:cover:cover}, notice that for every vertex $v \in V(G)$, there is some $C \in \mP_i$ such that $B_G(v, r/2 + 1) \subseteq C$. This implies that $v \in C'$, as desired. 
\end{proof}

    We record the following immediate consequence of Lemma \ref{lem: TFAE} for future reference: 
    
    \begin{corl}[Asymptotic dimension in terms of padded decompositions]\label{corl:pdtoBad}
        Let $G$ be a locally finite Borel graph and let $n \in \N$. Then $\Bad(G) \leq n$ if and only if for every large enough $r$, there exists some $\alpha > 1$ such that $G$ admits a Borel $(r, \alpha)$-padded decomposition with $n + 1$ layers.
    \end{corl}
    \begin{proof}
        Follows from Lemma~\ref{lem: TFAE} and Definition \ref{defn:Bad}.
    \end{proof}

    We conclude this subsection with the following definition, which will be used to simplify working with graphs of polynomial growth: 
    \begin{defn}[$(\bound, r_0)$-graphs]
        Let $\bound$, $r_0 > 0$. A graph $G$ is called a \emphdef{$(\bound, r_0)$-graph} if $\gamma_G(r) \leq r^\bound$ for all $r \geq r_0$, i.e., if $|B_G(v,r)| \leq r^\bound$ for all $r \geq r_0$ and $v \in V(G)$.
    \end{defn}
    
        Clearly, if $G$ is a $(\bound, r_0)$-graph, then $G$ is of polynomial growth and the asymptotic growth rate of $G$ is at most $\bound$. 
        Conversely, if $G$ is a graph with asymptotic growth rate $\bound$, then for every $\epsilon > 0$, there is some $r_0$ such that $G$ is a $(\bound + \epsilon, r_0)$-graph. 


\subsection{A randomized ball carving construction}\label{subsec:CSPsetup}
    Let $G$ be a Borel graph with finite maximum degree. To find a padded decomposition of $G$ with desired properties, we employ a randomized variant of a ``ball carving construction'' from computer science \ep{see, e.g., \cite{LinialSaks,Bartal, krauthgamer2007intrinsic,  Flitser,RG} for other applications of this technique}. To each positive integer $M$ and a function $t \colon V(G) \to \set{0,1, \ldots, M}$, 
we associate a $2M$-bounded partition $\mP_t$ of $V(G)$ as follows. Let $G^{2M}$ be the graph with vertex set $V(G)$ in which two vertices $x$, $y \in V(G)$ are adjacent if and only if $1 \leq \dist_G(x,y) \leq 2M$. This is a Borel graph of finite maximum degree, so, by Theorem~\ref{theo:Delta+1}, we may fix a Borel proper coloring of $G^{2M}$ with some finite number $k \in \N$ of colors. We interpret this coloring as a partition $V(G) = \bigcup_{i=0}^{k-1} I_i$, where each $I_i \subseteq V(G)$ is a Borel set such that $\dist_G(x, y) > 2M$ for all distinct $x$, $y \in I_i$. 
We process the sets $I_0$, \ldots, $I_{k-1}$ one by one and inductively define families $\C_0$, \ldots, $\C_{k-1}$ of finite subsets of $V(G)$ by:
\begin{equation}\label{eq:ballcarving}
    \C_0 \,\defeq\, \set{B_G(x, t(x)) \,:\, x \in I_0}, \quad \C_{i+1} \,\defeq\, \left\{B_G(x, t(x)) \backslash \left(\bigcup \C_0 \cup \ldots \cup \bigcup \C_i\right) \,:\, x \in I_{i+1}\right\}.
\end{equation}
Let $\mP_t \defeq \bigcup_{i=0}^{k-1} \C_i$. By construction, $\mP_t$ is a $2M$-bounded partition of $V(G)$ \ep{note that the sets in each $\C_i$ are pairwise disjoint because the points in $I_i$ are at distance greater than $2M$ from each other}. Furthermore, if the function $t \colon V(G) \to \set{0,1,\ldots, M}$ is Borel, then $\mP_t$ is a Borel partition.


    
    Our goal now is to argue, with the help of the \hyperref[thm: BLL]{Borel Local Lemma} (Theorem \ref{thm: BLL}) that for a certain choice of functions $t_1$, \ldots, $t_m$, the resulting partitions $(\mP_{t_1}, \ldots, \mP_{t_m})$ form an $(r, \alpha)$-padded decomposition of $G$ with $m$ layers for suitable parameters $r$, $\alpha$, and $m$. 

    We shall apply the \hyperref[thm: BLL]{Borel Local Lemma} in the generalized form described in Remark~\ref{remk:LLL}. That is, we will choose values $t(x) \in \set{0,1, \ldots, M}$ from a certain non-uniform probability distribution on $\set{0,1, \ldots, M}$. The distribution we shall use is given by 
\[\P[t = n] =\begin{cases}
 p(1-p)^n, &\textrm{if } n =0,\, 1,\, \ldots,\, M-1;\\
 (1-p)^M, &\textrm{if } n = M,
\end{cases}\]
where $p \in (0,1)$ is a parameter to be determined later. We call the resulting distribution the \emphd{truncated geometric distribution} and denote it by $\tgeo(p, M)$. The use of the truncated geometric distribution for the analysis of ball caving constructions was introduced in the paper \cite{LinialSaks} by Linial and Saks. It is well known that the geometric distribution is memoryless. For large $M$, the truncated distribution $\tgeo(p, M)$ is very close to the geometric distribution with parameter $p$, and as a result it is ``almost memoryless'': 

\begin{Lemma}
\label{lem:fact}
    Consider a random variable $t \sim \tgeo(p, M)$. We have:
    \begin{enumerate}[label=\ep{\normalfont\arabic*}]
        \item\label{item:fact:tail} For all nonnegative integers $n \leq M$, $\P[t \ge n] = (1 - p)^{n}$. 
        \item\label{item:fact:conditional} For all integers $m$, $n \ge 1$ such that $m + n < M$, $\P[t \le m + n \,|\, t \ge m] = 1 - (1 - p)^{n+1}$. 
    \end{enumerate}
\end{Lemma}
\begin{proof}
    \ref{item:fact:tail}
    $\P[t \ge n] = 1 - \sum_{i=0}^{n-1} p(1-p)^i = (1-p)^n$. 

\medskip
    
    \ref{item:fact:conditional} $\P[t \le m + n \,|\, t \ge m] = \dfrac{\P[m \le t \le m + n]}{\P[t \ge m]} = \dfrac{(1-p)^{m}-(1-p)^{m+n+1}}{(1-p)^{m}} = 1 - (1 - p)^{n+1}$.
\end{proof}

    Imagine running the ball carving construction described above with $t(x) \sim \tgeo(p, M)$ independently for all $x \in V(G)$.\footnote{The reader may feel somewhat uneasy about the idea of making an independent random choice for each $x \in V(G)$, since $G$ may have uncountably many vertices. It is useful to keep in mind here that in the statement of Theorem~\ref{thm: BLL}, probability is only involved in bounding $\P[A]$ for each constraint $A$, and this is a ``local'' calculation. In particular, the required probabilistic analysis can be performed on each connected component of $G$ separately.} We say that a ball $B_G(u, r)$ around a vertex $u \in V(G)$ is \emphdef{cut} if it intersects at least two different clusters of the resulting partition $\mP_t$.   
    The following lemma is the crucial probabilistic ingredient needed for the application of the \hyperref[thm: BLL]{Borel Local Lemma}:

\begin{Lemma}\label{lem:cut}
    Let $G$ be a $(\bound, r)$-graph with $\bound \geq 1$ and $r \geq 9$. 
    Let $u \in V(G)$ and suppose that 
    \[p \,\le\, \dfrac{1}{5 \bound} \quad \text{and} \quad M \,=\, \left\lfloor \dfrac{4\bound \log(\frac{1}{p})}{p} \right\rfloor.
    \]
    Then, with respect to the random choice of $t(x) \sim \tgeo(p, M)$ for all $x \in V(G)$, we have \[\P[\text{$B_G(u, r)$ is cut}] \,\le\, 20 r p.\]
\end{Lemma}

\begin{proof}
    Note that since $r \leq \ceil{r}$, $G$ is a $(\bound, \ceil{r})$-graph. Moreover, \[\P[\text{$B_G(u,r)$ is cut}] \,\leq\, \P[\text{$B_G(u,\ceil{r})$ is cut}].\] It follows that, since $\ceil{r} \leq r + 1 \leq 10r/9$ for $r \geq 9$, it is enough to prove the inequality
    \[
        \P[\text{$B_G(u, r)$ is cut}] \,\le\, 18 r p
    \]
    for integer $r$. Thus, from now on, we shall assume that $r$ is an integer. We may also assume that $r \leq 1/p$, since otherwise the conclusion of the lemma is trivial.
    
    For two vertices $x$, $y \in V(G)$ we write $y \prec x$ if $y \in I_i$ and $x \in I_j$ with $i < j$ (here $I_0$, \ldots, $I_{k-1}$ are the sets used in the inductive definition \eqref{eq:ballcarving}). 
    Set $B \defeq B_G(u, r)$ and let $B_x \defeq B_G(x, t(x))$ for all $x \in V(G)$.
    We say that the ball $B$ is \emphd{cut by $B_x$} if $B_y \cap B = \0$ for all $y \prec x$ and $\0 \neq B_x \cap B \neq B$. Clearly, $B$ is cut if and only if it is cut by some $B_x$.
    
    Define $A_{\mathrm{far}}$ to be the random event that there is some $B_x$ with $\dist_G(u, x) \ge M - r$ that cuts $B$. 
    Similarly, let $A_{\mathrm{near}}$ be the event that there is some $B_x$ with $\dist_G(u, x) < M - r$ that cuts $B$. 
    Then we have $\P[\text{$B$ is cut}]  \leq \P[A_{\mathrm{far}}] + \P[A_{\mathrm{near}}]$. We now bound $\P[A_{\mathrm{far}}]$ and $\P[A_{\mathrm{near}}]$ separately.

    To bound $\P[A_{\mathrm{far}}]$, we notice that for each $x \in V(G)$ with $\dist_G(u,x) \geq M - r$,
    \[
        \P[\text{$B_x$ cuts $B$}] \,\le\, \P[\dist_G(u, x) \le t(x) + r] \,\le\, \P[t(x) \ge  M - 2r] \,=\, (1 - p)^{M-2r}, 
    \]
    by Lemma~\ref{lem:fact}\ref{item:fact:tail}. Since we are assuming that $r \leq 1/p$, 
    \[
        (1 - p)^{M-2r} \,\le\, (1-p)^{-2/p} (1-p)^M.
    \]
    The function $z \mapsto (1-z)^{-2/z}$ is increasing for $0 < z < 1$, so, since $p \leq 1/5b \leq 1/5$, we obtain
    \[
        (1-p)^{-2/p} (1-p)^M \,\le\, (1-1/5)^{-10}(1-p)^{M} \,\leq\, 10 (1-p)^M.
    \]
    Using the inequality $1-z \leq e^{-z}$, valid for all $z \in \R$, and the bound $M \geq \dfrac{4\bound \log(\frac{1}{p})}{p} - 1$, we get
    \[
        10 (1-p)^M \,\le\, 10p^{4\bound} e^{p} \,\leq\, 15 p^{4\bound}.
    \]
    If $B_x$ cuts $B$, then $\dist_G(u,x) \leq M+r$, so the number of options for $x$ is at most 
    \begin{equation}\label{eq:LLlboundball}
        \gamma_G(M+r) \,\leq\, (M+r)^{\bound} \,\le\, \left (\frac{4\bound \log(\frac{1}{p})+1}{p} \right )^{\bound}.
    \end{equation}
    Observe that, since $5b \leq 1/p$ and $\log(1/p) \leq 1/p$, we have
    \[
        4\bound \log\left(\frac{1}{p}\right)+1 \,\leq\, 5\bound \log\left(\frac{1}{p}\right) \,\leq\, p^{-2}.
    \]
    Hence, the right-hand side of \eqref{eq:LLlboundball} is bounded above by $p^{-3b}$. 
    Finally, by the union bound, we have 
    \[\P[A_{\mathrm{far}}] \,\le\, p^{-3\bound} \, 15p^{4\bound} \,=\, 15p^\bound \,\le\, 15p. \]

    Now we turn to $\P[A_{\mathrm{near}}]$. 
    Observe that on the \ep{finite} set $B_G(u, M - r)$, the partial order $\prec$ is actually total, because all the vertices in $B_G(u, M-r)$ must belong to distinct sets $I_i$. 
    Let $X$ be the $\prec$-smallest element of $B_G(u, M - r)$ such that $B_X \cap B \neq \0$ (so $X$ is a random variable). Note that $X$ is well-defined since $B_u \cap B \supseteq \set{u} \neq \0$. By definition, for each $x \in B_G(u, M-r)$, $X = x$ if and only if the following two statements hold:
    \begin{itemize}
        \item $t(x) \geq \dist_G(u,x) - r$, and
        \item $t(y) < \dist_G(u, y) - r$ for all $y \in B_G(u, M-r)$ such that $y \prec x$.
    \end{itemize}
    If the event $A_{\mathrm{near}}$ takes place, then $B_X$ must cut $B$. It follows that if $t(X) \geq \dist_G(u,X) + r$, then $B_X \supseteq B$, and thus $A_{\mathrm{near}}$ does not happen. Therefore, for any $x \in B_G(u, M-r)$, we can write
    \begin{align*}
        \P[A_{\mathrm{near}}\,|\,X = x] \,&\le\,\, \P[t(x) < \dist_G(u, x) + r \,|\,  X = x] \\
        \,&=\, \P[t(x) < \dist_G(u, x) + r \,|\,  t(x) \ge  \dist_G(u, x) - r] \\
        &=\, 1-(1-p)^{2r} \\
        &\le\, 2rp,
    \end{align*}
    where the first equality holds since the values $t(y)$ for $y \neq x$ are independent from $t(x)$, and the second equality is satisfied by Lemma~\ref{lem:fact}\ref{item:fact:conditional}. 
    Therefore, we may conclude that 
    \begin{align*}
        \P[A_{\mathrm{near}}] 
        \,&=\, \sum_{x \in B_G(u, M-r)} \P[X = x]\;\P[A_{\mathrm{near}}\;|\;X = x] \\
        \,&\le\, \sum_{x \in B_G(u, M-r)} \P[X = x] \;2r p   \\
        \,&=\, 2rp.    
    \end{align*}
    
    Putting everything together, we obtain that
    \[\P[\text{$B$ is cut}]  \,\leq\, \P[A_{\mathrm{far}}] + \P[A_{\mathrm{near}}] \,\le\,  15p + 2rp \,\le\, 18rp,\]
    and the proof is complete.
\end{proof}

    We are now ready to define the CSP to which the \hyperref[thm: BLL]{Borel Local Lemma} will be applied. To this end, we observe that whether or not the ball $B_G(u,r)$ is cut by a partition $\mP_t$ is a property that is determined by the values of the function $t$ on the points in the finite set $B_G(u, M + r)$. Therefore, for each $m \in \N$, we can define a Borel CSP $\B_m \colon V(G) \to^? \set{0,1,\ldots,M}^m$ as follows: For each vertex $u \in V(G)$, let $A_{u,m}$ be the constraint with domain $\dom(A_{u,m}) = B_G(u, M+r)$ that is satisfied by a function $\bm{t} = (t_1, \ldots, t_m) \colon V(G) \to \set{0,1,\ldots,M}^m$ if and only if $B_G(u, r)$ is not cut in at least one of the partitions $\mathcal{P}_{t_1}$, \ldots, $\mathcal{P}_{t_m}$. Then we let\label{defn:theCSP}
    \[
        \B_m \,\defeq\, \set{A_{u,m} \,:\, u \in V(G)}.
    \]
    Following Remark \ref{remk:LLL}, we equip the set $\set{0,1,\ldots,M}^m$ with the product of $m$ copies of the distribution $\mathsf{tGeo}(p,M)$. That is, we imagine that for each $x \in V(G)$, the values $t_1(x)$, \ldots, $t_m(x)$ are chosen independently from the distribution $\mathsf{tGeo}(p,M)$. Then we have the following:
    
    \begin{Lemma}\label{lem:pd}
        Under the assumptions of Lemma~\ref{lem:cut}, \[\mathsf{p}(\B_m) \leq (20rp)^m \quad \text{and} \quad \mathsf{d}(\B_m) \leq  (2M+2r)^\bound - 1.\]
    \end{Lemma}
    \begin{proof}
        The bound on $\mathsf{p}(\B_m)$ follows by Lemma~\ref{lem:cut}, while the bound on $\mathsf{d}(\B_m)$ holds since
                \[
            \mathsf{d}(\B_m) \,\leq\, \gamma_G(2(M+r)) - 1 \,\leq\, (2M+2r)^\bound - 1,
        \]
        where we subtract $1$ because a constraint does not contribute to its own neighborhood in $\B_m$.
    \end{proof}

\subsection{Borel asymptotic dimension of graphs of polynomial growth}

After the preparations performed in \S\ref{subsec:CSPsetup}, we can now prove the following: 

\begin{theo}[Borel padded decomposition theorem I]
\label{thm:Borel_padded_decomposition_1}
    Let $G$ be a Borel $(\bound, r)$-graph with $\bound \geq 1$. Set $m \defeq \floor{\bound}+1$. Let $\epsilon > 0$ and define $\alpha \defeq (1+\epsilon)\frac{m}{m-\bound}$. If \[r \,>\, \max \{9, \, (1600 \bound^4)^{1/\alpha},\,  (8000\alpha \bound/\epsilon)^{2/\epsilon}\},\] then $G$ admits a Borel $(r, \alpha)$-padded decomposition with $m$ layers. 
\end{theo}

\begin{proof}
Consider the ball carving construction described in the previous subsection with parameters
\[
    p \,=\, \left(\dfrac{r^{\alpha}}{8 \alpha \bound \log r}\right)^{-1} \quad \text{and} \quad M \,=\, \left\lfloor\dfrac{4\bound \log(\frac{1}{p})}{p}\right\rfloor.
\]
Note that if $y \geq 2$ and $z \geq y^2$, then $z/\log z \geq y$. Applying this to $y = 40 b^2$ and $z = r^\alpha$ yields 
\[\frac{r^{\alpha}}{8 \alpha \bound \log r} \,\ge\, 5 \bound ,\] 
because $r > (1600 \bound^4)^{1/\alpha}$. 
Hence, $p \leq 1/5b$ and the assumptions of Lemma~\ref{lem:cut} are fulfilled. 

Let $\B_m$ be the Borel CSP defined in the previous subsection (on page \pageref{defn:theCSP}). Note that, since $G$ is of polynomial growth, the graph $G_{\B_m}$ associated to the CSP $\B_m$ is also of polynomial growth, and hence we may apply to it the \hyperref[thm: BLL]{Borel Local Lemma}. Specifically,
we observe that
\begin{equation}\label{eq:2M}
    2M \,\leq\, \dfrac{8\bound \log(\frac{1}{p})}{p} \,=\, \dfrac{r^\alpha \log(\frac{1}{p})}{\alpha \log r} \,\leq\, r^\alpha,
\end{equation}
where in the last step we use that $1/ p \leq r^\alpha$. Since $r \leq r^\alpha$ as well, Lemma~\ref{lem:pd} and Theorem~\ref{thm: BLL} imply that the CSP $\B_m$ has a Borel solution provided that
\[e \left(20rp\right)^{m}(2M+2r)^{\bound} \,<\, e \left(20rp\right)^{m}(5r^\alpha)^{\bound} \,<\, 1. \]
The last inequality is equivalent to 
\[\frac{r^{\alpha}}{8 \alpha \bound \log r} \,=\, 1/p \,>\, 20 \, e^{\frac{1}{m}}\,5^{\frac{\bound}{m}} \, r^{1+\frac{\alpha \bound}{m}}.\]
%
Using that $m > \bound \geq 1$, we obtain  $e^{1/m} 5^{\bound/m} \le 5^2 = 25$, so it suffices to prove that
\[\frac{r^{\alpha}}{8 \alpha \bound \log r} \,>\, 500 \, r^{1+\frac{\alpha \bound}{m}}.\]
Recall that $\alpha = (1+\epsilon)\frac{m}{m-\bound}$. Hence it is enough to get 
\[\frac{r^{\epsilon}}{\log r} \,>\, 4000 \alpha \bound.\] 
Notice that  $\log r \le \frac{2}{\epsilon}r^{\epsilon/2}$ (because $z > \log z$ for all $z > 0$), and hence it suffices to have \[r^{\epsilon/2} \,>\, \frac{8000\alpha\bound}{\epsilon},\] i.e., $r > (8000\alpha\bound/\epsilon)^{2/\epsilon}$, which holds by assumption. Therefore, the \hyperref[thm: BLL]{Borel LLL} may be applied, and we conclude that $\B_m$ has a Borel solution $\bm{t} = (t_1, \ldots, t_m) \colon V(G) \to \set{0,1,\ldots, M}^m$. We claim that the corresponding tuple $(\mP_{t_1}, \ldots, \mP_{t_m})$ of partitions of $G$ is a Borel $(r, \alpha)$-padded decomposition. Indeed, condition \ref{item:pd:padded} of Definition~\ref{def:padded_decomposition} holds since $\bm{t}$ is a solution to $\B_m$. By construction, each cluster in every $\mP_{t_i}$ has diameter at most $2M$. As $2M \leq r^\alpha$ by \eqref{eq:2M}, it follows that every $\mP_{t_i}$ is $r^\alpha$-bounded, and the proof is complete. 
\end{proof}

We are now ready to bound the Borel asymptotic dimension of graphs of polynomial growth:

\begin{corl}[Theorem~\ref{theo:Bad}]
    Every locally finite Borel graph $G$ satisfies $\Bad(G) \leq \ar(G)$.
\end{corl}

\begin{proof}
    We may assume that $\ar(G) < \infty$, since otherwise the statement is trivial. Also, if $\ar(G) < 1$, then in fact $\ar(G) = 0$ and all components of $G$ have uniformly bounded diameter, which implies $\Bad(G) = 0$. Thus, we may assume $\ar(G) \geq 1$. Let $m \defeq \lfloor \ar(G) \rfloor + 1$. Then $m > \ar(G)$, so we can pick some $\bound$ so that $\ar(G) < \bound < m$. By the definition of $\ar(G)$, 
    $G$ is an $(r, \bound)$-graph for all large enough $r$. 

Take any $\epsilon > 0$ and let $\alpha \defeq (1+\epsilon) m/ (m - b)$. By Theorem~\ref{thm:Borel_padded_decomposition_1}, for all sufficiently large $r$, $G$ admits a Borel $(r, \alpha)$-padded decomposition with $m$ layers. By Corollary~\ref{corl:pdtoBad}, this implies that \[\Bad(G) \,\leq\, m - 1 \,\leq\, \ar(G). \qedhere\] 
\end{proof}


    In the later sections we shall also need the following version of the Borel padded decomposition theorem, which allows us to construct $(r, \alpha)$-padded decompositions with $\alpha$ arbitrarily close to $1$, but with a larger number of layers: 

\begin{theo}[Borel padded decomposition theorem II]
\label{thm:Borel_padded_decomposition_2}

    Let $0< \epsilon \le \frac{1}{2}$, $\alpha \defeq 1 + \epsilon$. If $\bound \geq 1$ and
     \[r \,\geq\, (12000 \bound/\epsilon)^{4/\epsilon},\]
     then every Borel $(\bound, r)$-graph $G$ has a Borel $(r, \alpha)$-padded decomposition with $\ceil{\frac{6\bound}{\epsilon}}$ layers.
\end{theo}
\begin{proof}
    Note that $r > \max \set{9, (1600\bound^4)^{1/\alpha}}$ since $\alpha > 1$ and $\epsilon < 1$. Let $C \defeq 6/\epsilon$.  We use the same set-up as in the proof of Theorem~\ref{thm:Borel_padded_decomposition_1}, except with $m = \lceil Cb \rceil$. Following exactly the same analysis as in the proof of Theorem~\ref{thm:Borel_padded_decomposition_1}, the problem reduces to showing that
\[\frac{r^{\alpha}}{8 \alpha \bound \log r} \,>\, 20 \, e^{\frac{1}{C\bound}} \, 5^{\frac{1}{C}} \,r^{1+\frac{\alpha}{C}}.\]
    Using that $C\bound \ge \bound \geq 1$, we obtain  $e^{\frac{1}{C\bound}} \, 5^{\frac{1}{C}} \le 5^2 = 25$, so it suffices to prove that
\[\frac{r^{\alpha}}{8 \alpha \bound \log r} \,>\, 500 \, r^{1+\frac{\alpha}{C}}.\]
Recall that $C = 6/\epsilon$ and $\alpha = 1+\epsilon \leq 3/2$. Hence it is enough to get 
\[\frac{r^{\epsilon-\frac{1+\epsilon}{6/\epsilon}}}{ \log r} \,>\, 6000\bound. 
\]
Notice that  $\log r < \frac{2}{\epsilon}r^{\epsilon/2}$ (because $z > \log z$ for all $z > 0$), and hence it suffices to have 
\[r^{\epsilon/2-\frac{1+\epsilon}{6/\epsilon}} \,\geq\, \frac{12000 \bound}{\epsilon}. \]
When $0 < \epsilon \le 1/2$, we have 
$\epsilon/2 -\frac{1+\epsilon}{6 /\epsilon} = \epsilon / 3 - \epsilon^2 /6 \ge \epsilon / 3 - \epsilon /12  =  \epsilon/4$.
Hence, it is enough to have 
\[r^{\epsilon/4} \,\geq\, \frac{12000 \bound}{\epsilon}, \]
i.e.,
$r \geq (12000 \bound/\epsilon)^{4/\epsilon}$, which is true by assumption.
\end{proof}

\subsection{Strong padded decompositions}

    In the proof of Theorem \ref{theo:coarse}, we will require the following strengthening of the notion of a padded decomposition:

\begin{defn}[Strong padded decompositions]
\label{def:eta_strong_padded_decomposition}
    Let $G$ be a graph and fix parameters $0 < \eta < 1$ and $\alpha > 1$. 
    A tuple $(\mP_{1}, \mP_{2}, \ldots, \mP_{m})$ of $m$ partitions of $V(G)$ is called a \emphdef{$(1-\eta)$-strong $(r, \alpha)$-padded decomposition} of $G$ with $m$ \emphdef{layers} if the following properties are satisfied.
\begin{enumerate}[label=\ep{\normalfont\arabic*}]
    \item Each $\mP_i$ is $r^\alpha$-bounded. 
    \item For every $u \in V(G)$, we have
    \begin{equation}\label{eq:eta}
        |\set{i \in \set{1,\ldots, m} \,:\, \text{there is $C \in \mP_i$ such that $B_G(u,r) \subseteq C$}}| \,\geq\, (1-\eta)m.
    \end{equation}
\end{enumerate}
\end{defn}

    As the next lemma shows, given any Borel padded decomposition, we can modify it to form a Borel $(1-\eta)$-strong padded decomposition for any $\eta > 0$ by increasing the number of layers by a factor of $1/\eta$. The construction is essentially the same as in the proof of \cite[Lemma 2.1]{bowen2023definable} due to Bowen and Weilacher. Similar ideas appear in \cite[Lemma 1]{jorgensen2022geodesic} by J\o{}rgensen and Lang, but in their version of the result the right-hand side of \eqref{eq:eta} is just $2$ rather than $(1-\eta)m$. 
\begin{Lemma}[Strengthening a padded decomposition]
\label{lem: strengthen}
    Let $G$ be a locally finite Borel graph. Let $\alpha' > 1$, $r_0 > 0$, and $m \in \N^+$ be such that $G$ admits a Borel $(r, \alpha')$-padded decomposition with $m$ layers for every $r \geq r_0$. 
    Then for any $0 < \eta < 1$, $\alpha > \alpha'$, and
    \[
        r \,\ge\, \max\left\{r_0, \left(\frac{10m}{\eta}\right)^{\frac{\alpha'}{\alpha-\alpha'}}\right\},
    \]
    there exists a Borel $(1-\eta)$-strong $(r, \alpha)$-padded decomposition of $G$ with $\lceil\frac{m}{\eta}\rceil$ layers. 
\end{Lemma}
\begin{proof}
    Before starting the proof, let us review some notation. For a graph $G$ and a number $r \geq 1$, we let $G^r$ be the graph with vertex set $V(G)$ in which two vertices $u$, $v$ are adjacent if and only if $1 \leq \dist_G(u,v) \leq r$. For a subset $U \subseteq V(G)$, $G[U]$ is the subgraph of $G$ induced by $U$, i.e., the graph with vertex set $U$ and edge set $\set{uv \in E(G) \,:\, u, \, v \in U}$. We also let $U^\mathsf{c} \defeq V(G) \setminus U$ denote the complement of $U$.

    Now suppose the assumptions of the lemma hold. Let $N \defeq \lceil m/\eta\rceil - 1$. 
    By Lemma~\ref{lem: TFAE}\ref{item:from_pd_to_cover}, there exists a Borel $(4Nr, D)$-cover of $G$ with $m$ layers for $D \defeq (2Nr + 1)^{\alpha'}$, say $(\mU_1, \ldots, \mU_m)$. 
For each $1 \leq j \leq m$, define $U_j \defeq \bigcup \mU_j$. Let $H \defeq G^{2r}$. 
For all integers $0 \leq i \leq N$ and $1 \leq j \leq m$, let
\[
    S_i^j \,\defeq\, \{v \in V(G): \dist_H(v, U_j) = i\},
\]
and define $S_i \defeq \bigcup_{j = 1}^m S_i^j$. Let $\mathcal{F}_i$ be the set of all connected components of $H[S_i^\mathsf{c}]$.

    Take any $C \in \mathcal{F}_i$ and let $u \in C$. Since $U_1 \cup \ldots \cup U_m = V(G)$, there is some $1 \leq j \leq m$ such that $u \in U_j$. Let $C' \in \mU_j$ be the set containing $u$. We claim that $C \subseteq B_H(C', N)$. Indeed, otherwise there would be a vertex $v \in C$ with $\dist_H(v, C') > N$. Since $\dist_H(u, C') = 0$ and the graph $H[C]$ is connected, there must also exist a vertex $w \in C$ such that $\dist_H(w, C') = i$. Since $C \cap S_i = \0$, it follows that $\dist_H(w, U_j) \neq i$, and thus there is $x \in U_j \setminus C'$ such that $\dist_H(w, x) < i$. But then
    \[
        \dist_G(C', x) \,\leq\, 2r \,\dist_H(C', x) \,\leq\, 2r\,(\dist_H(C', w) + \dist_H(w,x)) \,<\, 4rN,
    \]
    contradicting the fact that the family $\mU_i$ is $4Nr$-disjoint. Now, since $C \subseteq B_H(C', N) = B_G(C', 2rN)$, it follows that the diameter of $C$ in $G$ is at most $D + 4rN$. We conclude that each $\mathcal{F}_i$ is a $(D + 4rN)$-bounded family of sets. Clearly, each $\mathcal{F}_i$ is also $2r$-disjoint (since there are no edges of $H$ between any two distinct sets in $\mathcal{F}_i$). And crucially, every vertex $u \in V(G)$ belongs to at least \[(N+1) - m \,\geq\, (1-\eta) (N+1)\] sets among $\mathcal{F}_0$, \ldots, $\mathcal{F}_N$, because for every $1 \leq j \leq m$, $u$ belongs to at most one set of the form $S_i^j$. 
    
    Now we apply Lemma~\ref{lem: TFAE}\ref{item:from_cover_to_pd} \ep{or, more accurately, the construction in the proof of that lemma} to $(\mathcal{F}_0, \ldots, \mathcal{F}_N)$ to obtain a $(1-\eta)$-strong $(r,\alpha)$-padded decomposition $(\mP_0, \ldots, \mP_N)$ of $G$ with $N+1$ layers, where $\alpha$ must satisfy $r^\alpha \geq D + 4rN + 2r$. After plugging in the value of $D$, this becomes
    \[
        r^\alpha \,\geq\, (2Nr + 1)^{\alpha'} + 4rN + 2r.
    \]
    For this to hold, it suffices to have $r^\alpha \geq (10Nr)^{\alpha'}$, which is true when $r \geq (\frac{10m}{\eta})^{\frac{\alpha'}{\alpha-\alpha'}}$.
\end{proof}

    Lemma~\ref{lem: strengthen} can be combined with Theorems~\ref{thm:Borel_padded_decomposition_1} and \ref{thm:Borel_padded_decomposition_2} to obtain strong padded decompositions of graphs of polynomial growth. We will specifically need the following result:

\begin{theo}[Borel strong padded decomposition theorem]
\label{thm:Borel_strong_padded_decomposition_2}
Fix $ 0 < \epsilon < 1$ and set $\alpha \defeq 1+ \epsilon$. If the parameters $0 < \eta < 1$, $r \geq 1$, $\bound \geq 1$, and $m \in \N^+$ satisfy 
\[r \,\geq\, \left(\frac{24000\bound}{\eta \epsilon}\right)^{\frac{8}{\epsilon}} \qquad \text{and} \qquad m \,=\, \left\lceil\frac{15\bound}{\eta\epsilon}\right\rceil,\]
then every Borel $(\bound, r)$-graph has a Borel $(1-\eta)$-strong $(r, \alpha)$-padded decomposition with $m$ layers. 
\end{theo}

\begin{proof}
    Let $G$ be a Borel $(\bound, r)$-graph. Applying Theorem~\ref{thm:Borel_padded_decomposition_2} with $\epsilon/2$ in place of $\epsilon$ shows that for \[R \,\geq\, r_0 \,\defeq\, (24000\bound/\epsilon)^{8/\epsilon},\] $G$ has a Borel $(R, 1 + \epsilon/2)$-padded decomposition with $\ceil{\frac{12\bound}{\epsilon}} \leq \frac{15\bound}{\epsilon}$ layers. Thus, by  Lemma~\ref{lem: strengthen}, 
for
\[
    r \,\geq\, \left(\frac{24000\bound}{\eta \epsilon}\right)^{\frac{8}{\epsilon}} \, \geq \, \max\left\{r_0,\, \left(\frac{150\bound}{\eta\epsilon}\right)^{\frac{3}{\epsilon}}\right\},
\]
there exists a Borel $(1-\eta)$-strong $(r, \alpha)$-padded decomposition with $\ceil{\frac{15\bound}{\eta\epsilon}}$ layers, as desired.
\end{proof}

 \section{Realizable cocycles}\label{sec:cocycle}
In this section we introduce the concept of a cocycle and establish a relationship between $\Z^n$-cocycles and maps to the graph $\shgr_{n, \infty}$. Throughout this section, $\G$ denotes a countable Abelian group, written additively. 

\begin{defn}[Cocycles]
Let $E$ be an equivalence relation on $X$ and let $\Gamma$ be an Abelian group. 
A $\Gamma$-\emphdef{cocycle} on $E$ is a function $\updelta \colon E \to \Gamma$ such that, for all $x \,E\, y \, E\, z$:
\begin{itemize}
    \item $\updelta (x, x) = 0$, and 
    \item $\updelta (x, y) + \updelta (y, z) = \updelta (x, z)$. 
\end{itemize}
We write $\updelta(u,v) = \infty$ if $(u,v) \not \in E$ (and when $\G = \Z^n$, we have $\|\updelta(u,v)\|_\infty = \infty$ as well in this case).
\end{defn}
    
    Note that the conditions defining a cocycle are linear. Hence, sums of cocycles are cocycles.

    \begin{defn}[Borel reductions between equivalence relations]
        Let $X$ and $Y$ be standard Borel spaces and let $E$, $F$ be equivalence relations on $X$, $Y$ respectively. A \emphdef{Borel reduction} from $E$ to $F$ is a Borel map $f \colon X \to Y$ such that for all $x$, $x' \in X$, we have
            \[x\,E\,x' \quad\iff\quad f(x)\,F\,f(x').\]
    \end{defn}
    
    Let $X$ be a standard Borel space and let $E$ be a Borel equivalence relation on $X$. Suppose that $\G$ is a countable Abelian group and $\G \acts Y$ is a free Borel action of $\G$ on a standard Borel space $Y$. Recall that the \emphd{orbit equivalence relation} $E(Y, \G)$ is the equivalence relation on $Y$ whose classes are the orbits of the action $\G \acts Y$. Given a Borel reduction $f \colon X \to Y$ from $E$ to $E(Y, \G)$, we define the \emphd{associated cocycle} $\updelta f \colon E \to \G$ on $E$ by making $\updelta f(x,x')$ be the unique group element $\gamma \in \G$ such that $\gamma \cdot f(x) = f(x')$. It is straightforward to check that $\updelta f$ is indeed a Borel $\G$-cocycle. We call cocycles that are constructed in this way \emph{realizable}:
    
    \begin{defn}[Realizable cocycles]\label{defn:realize}
        Let $E$ be a Borel equivalence relation on a standard Borel space $X$. A Borel $\G$-cocycle $\updelta \colon E \to \G$ is \emphd{realizable} if there exist a standard Borel space $Y$, a free Borel action $\G \acts Y$, and a Borel reduction $f \colon X \to Y$ from $E$ to $E(Y, \G)$ 
        such that $\updelta = \updelta f$.
    \end{defn}
    
    We will use realizable Borel cocycles to construct maps from Borel graphs to the graphs of the form $\shgr_{n, \infty}$ via the following lemma:
    
    \begin{Lemma}\label{lemma:map_to_free}
        Let $G$ be a locally finite Borel graph and suppose that $\updelta$ is a realizable $\Z^n$-cocycle on $\sim_G$. Then there exists a Borel map $f \colon V(G) \to \Free(2^{\Z^{n+1}})$ such that for all $u$, $v \in V(G)$,
        \[
            \dist_\infty(f(u), f(v)) \,=\, \|\updelta(u,v)\|_\infty.
        \]
    \end{Lemma}
    \begin{proof}
        Since the cocycle $\updelta$ is realizable, there exist a free Borel action $a \colon \Z^n \acts Y$ on a standard Borel space $Y$ and a Borel reduction $g \colon V(G) \to Y$ from $\sim_G$ to $E(Y, \Z^n)$ such that $\updelta = \updelta g$. In other words, for all $u$, $v \in V(G)$ such that $u \sim_G v$, we have
        \begin{equation}\label{eq:gcocycle}
            g(v) \,=\, \updelta(u,v) \cdot_a g(u).
        \end{equation}
        By Mycielski's theorem \cite[Theorem 19.9]{AnushDST}, there is an uncountable Borel set $T \subseteq \Free(2^\Z)$ that intersects every orbit of the shift action $\Z \acts \Free(2^\Z)$ in at most one point. By the Borel isomorphism theorem \cite[\S15]{KechrisDST}, we may assume, without loss of generality, that $Y = T$. Since $Y$ is now a subset of $\Free(2^\Z)$, there are two actions that can be applied to the elements of $Y$:
        \begin{itemize}
            \item the shift action $\Z \acts \Free(2^\Z)$, which we indicate by $\cdot$;
            
            \item the action $a \colon \Z^n \acts Y$, which we indicate by $\cdot_a$.
        \end{itemize}
        By the choice of $T = Y$, the following implication holds for all $y$, $y' \in Y$ and $i \in \Z$:
        \begin{equation}\label{eq:doubleshift}
            i \cdot y \,=\, y' \quad \Longrightarrow \quad i = 0 \text{ and } y = y'.
        \end{equation}
        Using the decomposition $\Z^{n+1} = \Z^n \times \Z$, we define a map $f \colon V(G) \to 2^{\Z^{n+1}}$ by 
        \begin{equation}\label{eq:f}
            f(u)(z, i) \, \defeq\, (z \cdot_a g(u))(i) \quad \text{for all } u \in V(G) \text{ and } z \in \Z^n, \ i \in \Z.
        \end{equation}
        Observe that for all $u$, $v \in V(G)$ and $z \in \Z^n$, $i \in \Z$,
        \begin{align*}
            (z,i) \cdot f(u) \,=\, f(v) \quad \iff \quad &\forall x \in \Z^n \, \forall j \in \Z, \ ((z,i) \cdot f(u))(x,j) \,=\, f(v)(x,j) \\
            \iff \quad &\forall x \in \Z^n \, \forall j \in \Z, \ f(u)(x+z,j + i) \,=\, f(v)(x,j) \\
            [\text{by \eqref{eq:f}}] \qquad\iff \quad &\forall x \in \Z^n \, \forall j \in \Z, \ ((x +z) \cdot_a g(u))(j + i) \,=\, (x \cdot_a g(v))(j) \\
            \iff \quad &\forall x \in \Z^n \, \forall j \in \Z, \ (i \cdot ((x + z) \cdot_a g(u)))(j) \,=\,  (x \cdot_a g(v))(j) \\
            \iff \quad &\forall x \in \Z^n, \ i \cdot ((x+z) \cdot_a g(u)) \,=\, x \cdot_a g(v) \\
            [\text{by \eqref{eq:doubleshift}}] \qquad \iff \quad & i = 0 \text{ and } z \cdot_a g(u) \,=\, g(v) \\
            [\text{by \eqref{eq:gcocycle}}]\qquad \iff \quad &i = 0 \text{ and } z = \updelta(u,v).
        \end{align*}
        This implies that $f(u) \in \Free(2^{\Z^{n+1}})$ for all $u \in V(G)$ \ep{by taking $u = v$}, that $f(u)$ and $f(v)$ are in the same component of the graph $\shgr_{n+1, \infty}$ if and only if $u \sim_G v$, and also that $\updelta f (u,v) = (\updelta (u,v), 0)$ for all $u \sim_G v$. In particular, for all $u \sim_G v$, we have
        \[
            \dist_\infty(f(u), f(v)) \,=\, \|(\updelta(u,v), 0)\|_\infty \,=\, \|\updelta(u,v)\|_\infty.
        \]
        Thus, the function $f$ has all the desired properties.
    \end{proof}
    

    Our goal now is to find necessary and sufficient conditions under which a given cocycle is realizable. To this end, we introduce some terminology. Let $X$ be a standard Borel space and $E$ be a Borel equivalence relation on $X$. We say that $E$ is \emphd{countable} \ep{resp.\ \emphd{finite}} if every $E$-class is countable \ep{resp.\ finite}. For an introduction to the theory of countable Borel equivalence relations, see the survey \cite{KechrisCBER} by Kechris. A Borel equivalence relation $E$ is \emphd{smooth} if it admits a Borel reduction to the equality relation on some standard Borel space. All finite Borel equivalence relations are smooth \cite[Lemma 5.21]{pikhurko2020borel}. A countable Borel equivalence relation is smooth if and only if it admits a Borel \emphd{transversal}, i.e., a Borel subset $T \subseteq X$ that contains exactly one point from every $E$-class. Other equivalent characterizations of smoothness for countable Borel equivalence relations can be found in \cite[\S3.5]{KechrisCBER}.
    
    

    Let $\Gamma$ be a countable Abelian group. We say that a $\G$-cocycle $\updelta \colon E \to \Gamma$ is \emphdef{smooth}\footnote{We remark that this notion is unrelated to the use of the term ``smooth cocycle'' in the paper \cite{miller2020existence} by Miller.} if the relation
\[
    E_\updelta \,\defeq\, \set{(x,x') \in E \,:\,\updelta(x,x') = 0}
\]
on $X$ is smooth. For example, this occurs when the relation $E_\updelta$ is finite. We can now state our characterization of realizable cocycles on countable Borel equivalence relations:




 
 
\begin{theo}\label{thm:cocycle-embedding}
    Let $E$ be a countable Borel equivalence relation on a standard Borel space $X$ and let $\Gamma$ be a countable Abelian group. A Borel $\G$-cocycle $\updelta$ on $E$ is realizable if and only if it is smooth. 
\end{theo}

\begin{proof}
    Assume first that $\updelta$ is realizable. Let $Y$ and $f$ be as in Definition \ref{defn:realize}, so $\updelta = \updelta f$. Then the map $f$ is a Borel reduction from $E_\updelta$ to the equality relation on $Y$, and therefore $E_\updelta$ is smooth.
    
    Now suppose that $E_\updelta$ is smooth. By \cite[Proposition 3.12]{KechrisCBER}, this implies that the quotient space $X' \defeq X/E_\updelta$ is standard Borel. For each $x \in X$, let $[x] \in X'$ denote the $E_\updelta$-class of $x$. Define an equivalence relation $E'$ on $X'$ by the formula
    \[
        [x] \, E'\, [x'] \quad\vcentcolon \Longleftrightarrow \quad x \, E \, x',
    \]
    and let $\updelta' \colon E' \to \Gamma$ be the cocycle given by $\updelta'([x], [x']) \defeq \updelta(x,x')$. Then $E_{\updelta'}$ is trivial (i.e., $[x] \, E_{\updelta'}\, [x']$ if and only if $[x] = [x']$), and since the quotient map $X \to X' \colon x \mapsto [x]$ is a Borel reduction from $E$ to $E'$, we may replace $X$, $E$, and $\updelta$ by $X'$, $E'$, and $\updelta'$ respectively and assume, without loss of generality, that $E_\updelta$ is trivial, i.e., $\updelta(x,x') = 0$ if and only if $x = x'$. 
    
    With this assumption, we define an equivalence relation $F$ on $\G \times X$ by
    \[
        (\gamma, x) \, F\, (\gamma', x') \quad \vcentcolon\Longleftrightarrow \quad \text{$x \, E \, x'$ and $\gamma - \gamma' = \updelta (x, x')$}.
    \]
    The fact that this is indeed an equivalence relation follows from the cocycle property of $\updelta$. Note that if $(\gamma, x) \,F\, (\gamma', x')$ and $\gamma = \gamma'$, then $x = x'$ as well, since in that case $\updelta (x, x') = \gamma - \gamma' = 0$ and we are assuming that $E_\updelta$ is trivial. In other words, the members of any $F$-equivalence class have distinct first coordinates. We now claim that $F$ is a smooth equivalence relation. List the elements of $\G$ in an arbitrary order as $\gamma_0$, $\gamma_1$, \ldots{} and define a preorder $\preceq$ on $\G \times X$ by
    \[(\gamma_i, x) \,\preceq\, (\gamma_j, x') \quad \vcentcolon\iff \quad i \,\le\, j.\]
    Then $\preceq$ induces a well-ordering on every $F$-class, so we can define a transversal $T$ for $F$ by picking the $\preceq$-smallest element in each class. A routine application of the \hyperref[theo:LN]{Luzin--Novikov theorem} \ep{Theorem~\ref{theo:LN}} shows that $T$ is Borel. 
    Therefore, $F$ is smooth and $Y \defeq (\G \times X)/F$ is a standard Borel space. 

    Let us now define a Borel action $\G \acts Y$ as follows:
    \[
        \delta \cdot [(\gamma, x)]_F \,\defeq\, [(\delta + \gamma, x)]_F.
    \]
    This is well-defined since if $(\gamma, x) \,F\, (\gamma', x')$, then \[(\delta + \gamma) - (\delta + \gamma') \,=\, \gamma -\gamma' \,=\, \updelta(x, x'),\] and hence $(\delta + \gamma, x) \,F\, (\delta + \gamma', x')$. Furthermore, this action is free since $(\gamma, x) \,F\, (\delta + \gamma, x)$ implies $\delta  = \updelta(x, x) = 0$.
    Now we define a function $f \colon X \to Y$ via
    $f(x) \defeq [(0, x)]_F$.
    This is a Borel reduction from $E$ to $E(Y, \G)$. Indeed, if $x \,E\, x'$, then
    \begin{equation}\label{eq:cocyclecompute}
        f(x') \,=\, \updelta(x,x') \cdot f(x),
    \end{equation}
    and hence $f(x)$ and $f(x')$ belong to the same $\G$-orbit. 
    Conversely, if $f(x)$ and $f(x')$ are in the same orbit, then there is a group element $\gamma \in \G$ with $f(x') = \gamma \cdot f(x)$, which implies that $(\gamma, x) \,F\, (0, x')$ and hence $x \, E\, x'$. Finally, \eqref{eq:cocyclecompute} implies that $\updelta f = \updelta$, and the proof is complete.
\end{proof}

\section{Coarse embeddings}\label{sec:coarse}

\subsection{Nested padded decompositions}

    In our proof of Theorem~\ref{theo:coarse} we apply Theorem~\ref{thm:Borel_strong_padded_decomposition_2} repeatedly for different values of $r$ in order to produce a family of padded decompositions that will be used to handle different distance scales in the construction. It will be important for us to ensure that these padded decompositions ``cohere'' with each other. In this subsection we formalize what that means and prove that such a ``coherent'' family of padded decompositions exists. 

    
    \begin{defn}[Refinements]
        Given two partitions $\mP$, $\mP'$ of a set $X$, we say that $\mP$ is a \emphdef{refinement} of $\mP'$, in symbols $\mP \preceq \mP'$, if for every $C \in \mP$, there is $C' \in \mP'$ such that $C \subseteq C'$.  
        
        Given two $m$-tuples $\D = (\mP_1, \ldots, \mP_m)$ and $\D' = (\mP'_1, \ldots, \mP'_m)$ of partitions of a set $X$, we say that $\D$ is a \emphdef{refinement} of $\D'$, in symbols $\D \preceq \D'$, if $\mP_i \preceq \mP'_i$ for all $i \in [m]$. 
    \end{defn}
    

\begin{Lemma}[Creating nested decompositions]
\label{lem:nesting}
    Let $G$ be a locally finite Borel graph. Let $r > 0$ and $0 < D \leq D'$ and suppose that $\D = (\mP_1, \ldots, \mP_m)$ and $\mF = (\mQ_1, \ldots, \mQ_m)$ are tuples of Borel partitions of $V(G)$ such that each $\mP_i$ is $D$-bounded and each $\mQ_i$ is $D'$-bounded.
    
    If $r \geq D$, then there exists a tuple $\D' = (\mP_1', \ldots, \mP_m')$ of $D'$-bounded Borel partitions of $V(G)$ such that $\D \preceq \D'$, and for any $v \in V(G)$ and $i \in [m]$, if there is a cluster $C \in \mQ_i$ with $B_G(v, 2r) \subseteq C$, then there also is a cluster $C' \in \mP'_i$ with $B_G(v, r) \subseteq C'$.
    
\end{Lemma}
\begin{proof}
    Consider the partitions $\mP_i$ and $\mQ_i$ for some $i \in [m]$. For a vertex $x \in V(G)$, let $C_x \in \mP_i$ be the cluster containing $x$. For each $C \in \mQ_i$, we define 
    \[
        C' \,\defeq\, \{x \in C \,:\, C_x \subseteq C\}.
    \]
    By construction, for any $x \in C'$, we have $C_x \subseteq C'$.
Let 
$V' \defeq \bigcup \{C' \,:\, C \in \mQ_i\}$.
Observe that if $x \in V(G) \setminus V'$, then $C_x \subseteq V(G) \setminus V'$, so the following is a partition of $V(G)$:
\[\mP'_i \,\defeq\, \{C' \,:\, C \in \mQ_i\} \cup \{C_x \,:\, x \in V(G) \setminus V'\}.\]
Clearly, $\mP'_i$ is Borel and $\mP_i \preceq \mP'_i$. Note that for each $C \in \mQ_i$ and $x \in V(G) \setminus V'$, \[\diam C' \,\leq\, \diam C \,\leq\, D' \qquad \text{and} \qquad \diam C_x \,\leq\, D \,\leq\, D'.\]
Thus, the partition $\mP'_i$ is $D'$-bounded. Now suppose that a vertex $v \in V(G)$ and a cluster $C \in \mQ_i$ satisfy $B_G(v, 2r) \subseteq C$. If $x \in B_G(v,r)$, then, since $\diam C_x \leq D \leq r$, we have $C_x \subseteq B_G(v, 2r) \subseteq C$, and so $x \in C'$. It follows that the tuple $\D' \defeq (\mP'_1, \ldots, \mP'_m)$ has all the desired properties.
%
%
\end{proof}


\begin{corl}[Nested decompositions in graphs of polynomial growth]
\label{corl: nested_decompositions}
Fix $0 < \epsilon < 1$ and set $\alpha \defeq 1+ \epsilon$. Suppose the parameters $0 < \eta < 1$, $r_0 \geq 1$, $\bound \geq 1$, and $m \in \N^+$ satisfy 
\[r_0 \,\geq\, \left(\frac{48000\bound}{\eta \epsilon}\right)^{\frac{16}{\epsilon}} \qquad \text{and} \qquad m \,=\, \left\lceil\frac{30\bound}{\eta\epsilon}\right\rceil,\]
and let $(r_n)_{n \in \N}$ be a sequence of integers such that $r_{n+1} \geq r_n^\alpha$ for all $n \in \N$. Then for every Borel $(\bound, r_0)$-graph $G$, there is a sequence $(\D_n)_{n \in \N}$ in which every $\D_n$ is a Borel $(1-\eta)$-strong $(r_n, \alpha)$-padded decomposition with $m$ layers and $\D_n \preceq \D_{n+1}$. 
%
%
\end{corl}

\begin{proof}
    Set $\alpha' \defeq 1 + \epsilon/2$ and define $D_n \defeq (2r_n)^{\alpha'}$ for all $n \in \N$. Observe that
    \[D_{n} \,=\, (2r_n)^{\alpha'} \,\leq\, 4 r_n^{\alpha'} \,\leq\, r_n^\alpha \,\leq\, r_{n+1},\] since $r_n \geq r_0 \geq 4^{2/\epsilon}$. We inductively construct a sequence $(\D_n)_{n \in \N}$, where each $\D_n = (\mP_{n,1}, \ldots, \mP_{n,m})$ is a tuple of $D_n$-bounded Borel partitions of $G$ and $\D_n \preceq \D_{n+1}$ for all $n \in \N$. To begin with, by Theorem~\ref{thm:Borel_strong_padded_decomposition_2}, there exists a Borel $(1-\eta)$-strong $(r_0, \alpha')$-padded decomposition $\D_0$ of $G$ with $m$ layers. Note that each layer of $\D_0$ is $D_0$-bounded since $D_0 = (2r_0)^{\alpha'} > r_0^{\alpha'}$.
    
    Once $\D_n$ has been defined, we construct $\D_{n+1}$ as follows. By Theorem~\ref{thm:Borel_strong_padded_decomposition_2}, there exists a Borel $(1-\eta)$-strong $(2r_{n+1}, \alpha')$-padded decomposition $\mF_{n+1} = (\mQ_{n+1, 1}, \ldots, \mQ_{n+1, m})$ of $G$. Notice that every layer of $\mF_{n+1}$ is $D_{n+1}$-bounded. Since $r_{n+1} \geq D_n$, we may apply Lemma~\ref{lem:nesting} with $r = r_{n+1}$, $D = D_n$, and $D' = D_{n+1}$ to get a tuple $\D_{n+1} = (\mP_{n+1,1}, \ldots, \mP_{n+1,m})$ of Borel partitions such that
    \begin{itemize}
        \item $\D_n \preceq \D_{n+1}$,
        \item each $\mP_{n+1,i}$ is $D_{n+1}$-bounded \ep{and hence $r_{n+1}^\alpha$-bounded}, and
        \item for any $v \in V(G)$ and $i \in [m]$, if there exists a cluster $C \in \mQ_{n+1,i}$ with $B_G(v, 2r_{n+1}) \subseteq C$, then there also exists a cluster $C' \in \mP_{n+1,i}$ with $B_G(v, r_{n+1}) \subseteq C'$.
    \end{itemize}
    It follows that $\D_{n+1}$ is a $(1-\eta)$-strong $(r_{n+1}, \alpha)$-padded decomposition, as desired.
    %
\end{proof}

\subsection{Constructing coarse embeddings into grids}

    We have now established all the preliminaries necessary to begin the proof of Theorem~\ref{theo:coarse}. We start with the following technical definition: 
    \begin{defn}[Dumpling-contractions]
        For a graph $G$, a map $\phi \colon V(G) \to \Z^n$ is called a \emphdef{contraction} if $\norm{\phi(u) - \phi(v)}_{\infty} \le 1$ for every edge $uv \in E(G)$. 
        For a set $C \subseteq V(G)$, an \emphd{$n$-dimensional dumpling-contraction} for $C$, or \emphdef{$\dumpling$-contraction} for short, is a contraction $\phi \colon V(G) \to \Z^n$ such that $\phi(u) = 0$ for all $u \in \partial C$.
        For a partition $\mP$ of $G$, a map $\phi \colon V(G) \to \Z^n$ is a \emphdef{$\dumpling$-contraction over $\mP$} if it is a $\dumpling$-contraction for every cluster $C \in \mP$.  
        Finally, for a tuple $\D = (\mP_1, \mP_2, \ldots, \mP_m)$ of partitions, a map $\phi \colon V(G) \to \Z^n$ is a \emphdef{$\dumpling$-contraction over $\D$} if it can be expressed as $\phi = \bigoplus_{i = 1}^m \phi_i$, where for each $i \in [m]$, $\phi_i \colon V(G) \to \Z^{n_i}$ is a $\dumpling$-contraction over $\mP_i$.
    \end{defn}

    \begin{remk*}
        We use the term ``dumpling'' to indicate that the boundary of a set is contracted to a single point, in analogy with the shaping process of a dumpling, as represented by the 
        symbol $\dumpling$. 
    \end{remk*}

    To facilitate our inductive construction, we will need to take quotients of graphs with respect to their partitions:

\begin{defn}[Quotient  graphs]
Let $G$ be a graph and let $\mP$ be a partition of $V(G)$. The \emphdef{quotient graph} of $G$ with respect to $\mP$ is the graph $G/\mP$ with vertex set $\mP$ and adjacency relation
\[
    \{(C, C') \in \mP^2 \,:\, C \neq C' \text{ and } E(G) \cap (C \times C') \neq \0\}.
\]
\end{defn}

    Note that if $G$ is a Borel graph and $\mP$ is a Borel partition of $V(G)$ into finite clusters, then the quotient graph of $G$ with respect to $\mP$ is also a Borel graph. 

    We now state the main technical lemma which will be used iteratively to prove Theorem~\ref{theo:coarse}: 

\begin{Lemma}
\label{lem:main}
    Fix real numbers $\alpha$, $\beta$, $\gamma$, $\epsilon$, $\bound$, and $r$ satisfying the following:
    \begin{enumerate}[label=\ep{\normalfont\arabic*}]
        \item\label{item:order} $\epsilon < 1 < \alpha < \beta < \gamma$ and $\bound$, $r \geq 1$;
        \item $r^{\gamma \epsilon/2} \ge 64$;
        \item $r^{(\gamma -\beta) \bound} \geq 2$.
    \end{enumerate}
    Let $\eta \in (0, \frac{1}{2})$ and let $G$ be a Borel $(\bound,r)$-graph. Let 
    $\D = (\mP_{1},\ldots, \mP_{m})$ and $\mF = (\mQ_{1}, \ldots, \mQ_{m})$ be Borel $m$-tuples of partitions of $G$ such that:
\begin{enumerate}[label=\ep{\normalfont\arabic*},resume]
    \item\label{item:LLL} $(1-\eta)m (\beta/\alpha - 1 - \gamma(1-\epsilon/2)) \geq 6\gamma \bound$;
    \item\label{item:DlessthanF} $\D \preceq \mF$;
    \item\label{item:Pbounded}  every $\mP_i$ is $r$-bounded; 
    \item\label{item:Qbounded} every $\mQ_i$ is $r^\beta$-bounded; 
    \item\label{item:eta} for all $u \in V(G)$, $B_G(u, r^{\beta/\alpha})$ is contained in a $\mQ_i$-cluster for at least $(1-\eta) m$ values of $i$.
\end{enumerate}
Let $\phi = \bigoplus_{i=1}^m \phi_i$ be a Borel function, where $\phi_i \colon V(G) \to \Z$ for all $i \in [m]$. 
Then there exists 
a Borel $\dumpling$-contraction $\psi = \bigoplus_{i=1}^m \psi_i$ over $\mF$, where $\psi_i \colon V(G) \to \Z$ for all $i \in [m]$, 
such that:
\begin{enumerate}[label=\ep{\itshape\alph*}]
    \item\label{item:a}  $\psi_{i}$ is constant on $\mP_i$-clusters for all $i \in [m]$;
    \item\label{item:b} for all $u$, $v \in V(G)$ with $r^{\beta} < \dist_G(u, v) \le r^{\gamma}$, the inequality
    \[
        |(\phi_i + \psi_i)(u) - (\phi_i + \psi_i)(v)| \,\ge\, r^{\gamma(1 - \epsilon)}
    \]
    is satisfied for at least $(1-\eta)m/2$ values of $i \in [m]$.
\end{enumerate}
\end{Lemma}

\begin{proof}
    When reading the following argument, it may be helpful to think of $\alpha$ as very close to $1$, while $\beta$ and $\gamma$ as large and close to each other, as that will be the case when we apply Lemma~\ref{lem:main} in the proof of Theorem~\ref{theo:coarse} (or, more precisely, Theorem~\ref{theo:coarse_p}) later on.

For $i \in [m]$, let $G_i \defeq G/\mP_i$. 
Note that $V(G_i) = \mP_i \preceq  \mQ_i$. Hence, we have a partition \[\mQ_i/\mP_i \,\defeq\, \{S/\mP_i \,:\, S \in \mQ_i\}\] of $\mP_i$, where $S/\mP_i = \{C \in \mP_i \,:\, C \subseteq S\}$ is the induced partition of $S$. For a vertex $u \in V(G)$, we let $[u]_{\mP_i}$ \ep{resp.\ $[u]_{\mQ_i}$} denote the cluster of the partition $\mP_i$ \ep{resp.\ $\mQ_i$} containing $u$. Since $\mP_i$ is $r$-bounded, for every two vertices $u$, $v \in V(G)$, we have
\begin{equation}\label{eq:hop}
    \dist_G(u, [v]_{\mP_i})\,\leq\, (r+1) \, \dist_{G_i}([u]_{\mP_i}, [v]_{\mP_i}).
\end{equation}
For a cluster $S \in \mQ_i$, let $\partial_i(S)$ denote the boundary of $S/\mP_i$ in the graph $G_i$.

    Let $I \defeq \set{k \in \Z \,:\, 0 \leq k \leq r^{\beta/\alpha}/(r+1)}$. 
    Given a function $t \colon \bigsqcup_{i = 1}^m \mQ_i \to I$, we construct a mapping $\psi^t \colon V(G) \to \Z^m$ in the following way. For $S \in \mQ_i$, let 
    \[
        \partial^t_i(S) \,\defeq\, \{C \in S/\mP_i \,:\, \dist_{G_i}(C, \partial_i(S)) \le t(S)\}.
    \]
    Given a vertex $u \in V(G)$, we define
    \[
        \psi_i^t (u) \,\defeq\, \dist_{G_i}([u]_{\mP_i}, \, \partial_i^t([u]_{\mQ_i})).
    \]
    Equivalently, we can write
    \[
        \psi_{i}^t(u) \,=\, \max \set{0, \, \dist_{G_i}([u]_{\mP_i}, \,\partial_i([u]_{\mQ_i})) - t([u]_{\mQ_i})}.
    \]
    %
Set $\psi^t \defeq \psi_1^t \oplus \cdots \oplus \psi_m^t$. Note that if the function $t$ is Borel, then $\psi^t$ is also Borel. Clearly, regardless of the choice of $t$, $\psi^t$ is constant on $\mP_i$-clusters, and it is a $\dumpling$-contraction over $\mQ_i$. Thus, condition \ref{item:a} is satisfied by $\psi^t$ for all $t$. It remains to argue that for some Borel function $t$, the corresponding mapping $\psi^t$ satisfies \ref{item:b}. To this end, we shall employ the \hyperref[thm: BLL]{Borel Local Lemma}.

    For any $u \in V(G)$, the value $\psi^t(u)$ is determined by the values $t([u]_{\mQ_i})$ for $i \in [m]$. Therefore, we can define a Borel CSP $\B : \bigsqcup_{i = 1}^m \mQ_i \to^{?} I$ as follows: For each pair of vertices $u$, $v \in V(G)$ such that $r^{\beta} < \dist_G(u, v) \le r^{\gamma}$, let $A_{u, v}$ be the constraint with domain $\{[u]_{\mQ_i}, [v]_{\mQ_i}: i \in [m]\}$ that is satisfied by a function $t: \bigsqcup_{i = 1}^m \mQ_i \to I$ if and only if there exist at least $(1-\eta)m/2$ indices $i$ 
such that
 \begin{equation}\label{eq:difference}
    |(\phi_{i}(u) + \psi_i^t (u)) - (\phi_{i}(v) + \psi_i^t(v))| \,\ge\, r^{\gamma(1 - \epsilon)}.
\end{equation}
By construction, if $t$ is a solution to $\B$, then the function $\psi^t$ has all the desired properties. Since the graph $G$ is of polynomial \ep{hence subexponential} growth, the graph $G_{\B}$ associated to the CSP $\B$ is also of subexponential growth, and thus we may apply the \hyperref[thm: BLL]{Borel Local Lemma} to it. Hence, to complete the proof, we need to verify the inequality
\[
    e \cdot \p(\B) \cdot (\d(\B) + 1) \,<\, 1.
\]

To bound $\d(\B)$, consider two distinct constraints $A_{u,v}$, $A_{u',v'} \in \B$ such that \[\dom(A_{u,v}) \cap \dom(A_{u', v'}) \,\neq\, \0.\] Then for some $i \in [m]$, either $u'$ or $v'$, for concreteness $u'$, belongs to the same $\mQ_i$-cluster as either $u$ or $v$. Since $G$ is a $(\bound, r)$-graph, each $\mQ_i$ is $r^\beta$-bounded, and $\dist_G(u', v') \leq r^\gamma$, this implies that
\[
    |N(A_{u,v})| \,\leq\, 2 r^{\beta \bound} r^{\gamma \bound} - 1 \,\leq\, r^{2\gamma \bound} - 1,
\]
where we use that $r^{(\gamma - \beta)\bound} \geq 2$ and subtract $1$ to exclude the case $A_{u,v} = A_{u',v'}$. Therefore,
\[
    \d(\B) \,\leq\, r^{2\gamma \bound} - 1.
\]

Now we turn our attention to $\p(\B)$. Suppose that the values $t(S)$ for $S \in \bigsqcup_{i=1}^m \mQ_i$ are chosen independently and uniformly at random from the set $I$. Consider any two vertices $u$, $v \in V(G)$ such that $r^{\beta} < \dist_G(u, v) \le r^{\gamma}$.
Since $\dist_G(u,v) > r^\beta$, condition \ref{item:Qbounded} implies that for every $i \in [m]$, $u$ and $v$ belong to different $\mQ_i$-clusters, and hence $[u]_{\mP_i}$ and $[v]_{\mP_i}$ reside in different clusters of $\mQ_i/\mP_i$. 
By~\ref{item:eta}, the ball $B_G(u, r^{\beta/\alpha})$ is contained in $[u]_{\mQ_i}$ for at least $(1-\eta)m$ values of $i$. 
By \eqref{eq:hop}, for any such $i$, $B_{G_i}([u]_{\mP_i}, r^{\beta/\alpha}/(r+1))$ is contained in $[u]_{\mQ_i}/\mP_i$, i.e.,
\[
    \dist_{G_i}([u]_{\mP_i}, \, \partial_i([u]_{\mQ_i})) \,\geq\, \left\lfloor \frac{r^{\beta/\alpha}}{r+1} \right\rfloor \,=\, \max I.
\]
If follows that the value $(\phi_{i} + \psi^t_{i})(u)$ is chosen uniformly at random from a set of $|I|$ distinct values, 
independently of $(\phi_{i} + \psi^t_{i})(v)$. Therefore,
\[
    \P[|(\phi_{i} + \psi_{i}^t)(u) - (\phi_{i} + \psi_{i}^t)(v)| < r^{\gamma(1-\epsilon)}] \,\le\, \frac{2r^{\gamma(1-\epsilon)}}{|I|}
\,\le\, \frac{2r^{\gamma(1-\epsilon)}}{r^{\beta/\alpha}/(r+1)}.
\]
Using that $r+1 \leq 2r$, 
we write
\[
\frac{2r^{\gamma(1-\epsilon)}}{r^{\beta/\alpha}/(r+1)}\,\leq\, \frac{4r^{\gamma(1-\epsilon)}}{r^{\beta/\alpha - 1}}. 
\]
Since the values $t(S)$ for clusters $S$ from each of the partitions $\mQ_1$, \ldots, $\mQ_m$ are chosen independently from each other, we conclude that the probability that \eqref{eq:difference} fails for a specific set of $\ceil{(1-\eta)m/2}$ indices $i$ such that $B_G(u, r^{\beta/\alpha}) \subseteq [u]_{\mQ_i}$ is at most
\[
    \left(\frac{4r^{\gamma(1-\epsilon)}}{r^{\beta/\alpha - 1}}\right)^{(1-\eta)m/2}.
\]
Taking the union over all such sets of indices, we obtain the bound
\[
    \P[A_{u,v}] \,\leq\, 2^m \, \left(\frac{4r^{\gamma(1-\epsilon)}}{r^{\beta/\alpha - 1}}\right)^{(1-\eta)m/2} \,\leq\, \left(\frac{64r^{\gamma(1-\epsilon)}}{r^{\beta/\alpha - 1}}\right)^{(1-\eta)m/2} \,\leq\, \left(\frac{r^{\gamma(1-\epsilon/2)}}{r^{\beta/\alpha - 1}}\right)^{(1-\eta)m/2} \,\leq\, r^{-3\gamma \bound},
\]
where we use that $\eta < 1/2$, $r^{\gamma \epsilon/2} \ge 64$, and the inequality  \ref{item:LLL}.
Therefore,
\[
    \p(\B) \,\leq\, r^{-3\gamma \bound}.
\]

From the above calculations, it follows that 
$\B$ has a Borel solution as long as 
\[e \cdot \p(\B) \cdot (\d(\B)+1) \,\le\, e \, r^{-3\gamma\bound} \, r^{2\gamma \bound} \,=\, e \, r^{-\gamma \bound}
\,<\, 1.\]
This inequality holds since $r^{\gamma \bound} \geq r^{\gamma \epsilon} \geq 4096$, and the proof is complete. 
\end{proof}

Finally, we can prove the following numerically explicit version of Theorem~\ref{theo:coarse}:

 \begin{theo}[Explicit form of Theorem~\ref{theo:coarse}]\label{theo:coarse_p}
        Let $G$ be a Borel $(\bound, r)$-graph with $\bound$, $r \geq 1$. Then for every $0 < \epsilon < \frac{1}{2}$, there exists a Borel function $f \colon V(G) \to \Free(2^{\Z^N})$, where 
        \[N \,=\, \left\lceil \frac{10^7 \, \bound \log(1/\epsilon)}{\epsilon^2} \right\rceil, 
        \]
        such that for all $u$, $v \in V(G)$, $\dist_\infty(f(u), f(v)) \leq \dist_G(u,v)$ and if \[\dist_G(u,v) \,\geq\, \max \left \{ \left(\frac{10^7 \, \bound}{\epsilon}\right)^{\frac{3000}{\epsilon^2}}, \, r^{\frac{15}{\epsilon}}\right\},
        \]
        then $\dist_\infty(f(u), f(v)) \geq \dist_G(u, v)^{1-\epsilon}$.
    \end{theo}

\begin{proof}
%
We introduce the following parameters:
\begin{align*}
    \alpha \,\defeq \, 1+\frac{\epsilon}{12}, \quad \beta \,\defeq\, \frac{12}{\epsilon}, \quad \gamma \,\defeq\, \alpha \beta \,=\, 1+\frac{12}{\epsilon}, \quad m \,\defeq\, \left\lceil \frac{1440 \bound}{\epsilon}\right\rceil.
\end{align*}
Let $x$ be an arbitrary number such that
\[
    x \,\geq\, \max \left\{\left(\frac{10^7 \, \bound}{\epsilon}\right)^{\frac{200}{\epsilon}}, \, r\right\}.
\]
Define a sequence of numbers $r_n$ as follows:
\[
    r_0 \,\defeq\, x, \quad r_{n+1} \,\defeq\, r_n^{\beta} \text{ for all } n \in \N \quad \text{(so $r_n = x^{\beta^n}$)}.
\]
Let $R_n \defeq r_{2n + 1}^{\alpha/\beta} = x^{\alpha \beta^{2n}}$ and define $\interval_n \defeq \left(R_n^\beta,\, R_n^\gamma\right]$. 
Since $\beta > \alpha$, Corollary~\ref{corl: nested_decompositions} applied with $\epsilon/12$ in place of $\epsilon$ and with $\eta = 1/4$ yields a sequence $(\D_n, \mF_n)_{n \in \N}$, where for all $n \in \N$,
\begin{itemize}
    \item $\D_n$ is a Borel $3/4$-strong $(r_{2n}, \alpha)$-padded decomposition of $G$ with $m$ layers,
    \item $\mF_n$ is a Borel $3/4$-strong $(r_{2n+1}, \alpha)$-padded decomposition of $G$ with $m$ layers, and
    \item $\D_n \preceq \mF_{n} \preceq \D_{n+1}$.
\end{itemize}
Write $\D_n = (\mP_{n,1}, \ldots, \mP_{n,m})$ and $\mF_n = (\mQ_{n,1}, \ldots, \mQ_{n,m})$. Now we inductively define a sequence of Borel functions $\phi_n \colon V(G) \to \Z^m$ as follows. Start with $\phi_0 \defeq 0$. Once $\phi_n = \bigoplus_{i=1}^m \phi_{n,i}$ has been defined, where $\phi_{n,i} \colon V(G) \to \Z$ for all $i \in [m]$, we wish to apply Lemma~\ref{lem:main} with
$R_n$, $\D_n$, $\mF_{n}$, $1/4$, and $\phi_n$ in place of $r$, $\D$, $\mF$, $\eta$, and $\phi$ respectively.
To this end, we observe that:

\smallskip

\begin{enumerate}[label=\ep{\normalfont\arabic*}]
    \item $\epsilon < 1 < \alpha < \beta < \gamma$ and $\bound$, $R_n \geq 1$;

    
    \item $R_n^{\gamma \epsilon/2} = R_n^{6 + \epsilon/2} \geq x > 64$;

    
    \item $R_n^{(\gamma - \beta) \bound} = R_n^\bound \geq x > 2$;

    
    \item $(3/4)m (\beta/\alpha - 1 - \gamma(1-\epsilon/2)) = (3/4)m (3 + \epsilon/2 + \epsilon/(\epsilon + 12)) > 9m/4 > 6\gamma \bound$;

    
    \item $\D_n \preceq \mF_{n}$;

    
    \item every cluster in $\D_n$ is $r_{2n}^{\alpha}$-bounded, and $r_{2n}^{\alpha} = r_{2n+1}^{\alpha/\beta} = R_n$;
    

    \item every cluster in $\mF_{n}$ is $r_{2n+1}^{\alpha}$-bounded, and $r_{2n+1}^\alpha = R_n^\beta$;
    
    
    \item for $u \in V(G)$, $B_G(u, r_{2n+1}) = B_G(u, R_n^{\beta/\alpha})$ is in a cluster of at least $3m/4$ partitions in $\mF_{n}$.
\end{enumerate}

\smallskip

\noindent Therefore, the assumptions of Lemma~\ref{lem:main} are satisfied, and hence there exists a Borel $\dumpling$-contraction $\psi_n = \bigoplus_{i=1}^m \psi_{n,i}$ over $\mF_{n}$, where $\psi_{n,i} \colon V(G) \to \Z$ for all $i \in [m]$, 
such that:


\begin{enumerate}[label=\ep{\itshape\alph*}]
    \item\label{item:a1}  $\psi_{n,i}$ is constant on $\mP_{n,i}$-clusters for all $i \in [m]$;
    \item\label{item:b1} for all $u$, $v \in V(G)$ with $\dist_G(u, v) \in \interval_n$, the inequality
    \[
        |(\phi_{n,i} + \psi_{n,i})(u) - (\phi_{n,i} + \psi_{n,i})(v)| \,\ge\, R_n^{\gamma(1 - \epsilon)} \,\geq\, \dist_G(u, v)^{1-\epsilon}
    \]
    is satisfied for at least $3m/8$ values of $i \in [m]$.
\end{enumerate}


\noindent We then let $\phi_{n+1} \defeq \phi_n + \psi_n$. In other words, we have $\phi_{n+1} = \sum_{k = 0}^n \psi_k$.

Next we use the techniques developed in \S\ref{sec:cocycle} to take the ``limit'' of the sequence $(\phi_n)_{n \in \N}$. For every pair of vertices $u$, $v \in V(G)$ with $u \sim_G v$, let
\[
    \updelta_n(u, v) \,\defeq\, \psi_n(v) - \psi_n(u).
\]
Then $\updelta_n$ is a Borel $\Z^m$-cocycle on $\sim_G$. Note that
\begin{equation}\label{eq:edge_bound1}
    \|\updelta_n(u,v)\|_\infty \,=\, \|\psi_n(v) - \psi_n(u)\|_\infty \,\leq\, 1 \quad \text{for all } uv \in E(G),
\end{equation}
because $\psi_n$ is a contraction. Let $\updelta_{n,i} \colon {\sim_G} \to \Z$ denote the $i$-th coordinate of $\updelta_n$. We claim that if $uv \in E(G)$ is an edge such that $\updelta_{n,i}(u,v) \neq 0$, then $u$ and $v$ belong to different $\mP_{n,i}$-clusters but the same $\mQ_{n, i}$-cluster. Indeed, if $u$ and $v$ are in the same $\mP_{n,i}$-cluster, then $\updelta_{n,i}(u,v) = 0$ as $\psi_{n,i}$ is constant on the $\mP_{n,i}$-clusters. On the other hand, if $u$ and $v$ are in different $\mQ_{n, i}$-clusters, then both $u$ and $v$ belong to the boundaries of their corresponding $\mQ_{n, i}$-clusters. It follows that $\updelta_{n, i}(u,v) = 0$ as $\psi_{n,i}$ is a $\dumpling$-contraction over $\mQ_{n, i}$. Since $\mP_{n,i} \preceq \mQ_{n,i} \preceq \mP_{n+1,i}$ for all $n \in \N$, there can be at most one index $n$ such that $u$ and $v$ belong to different $\mP_{n,i}$-clusters but the same $\mQ_{n, i}$-cluster. Therefore,
\begin{equation}\label{eq:edge_bound}
    |\set{n \in \N \,:\, \updelta_{n,i}(u,v) \neq 0}| \,\leq\, 1 \quad \text{for all } uv \in E(G).
\end{equation}
For any two vertices $u \sim_G v$, there is a finite $uv$-path in $G$. Since $\updelta_n(u,v)$ is equal to the sum of the values of $\updelta_n$ on the edges of this path, we have
\[
    |\set{n \in \N \,:\, \updelta_{n,i}(u,v) \neq 0}| \,<\,\infty \quad \text{for all } uv \in E(G).
\]
It follows that we can define a $\Z^m$-cocycle $\updelta \colon {\sim_G} \to \Z^m$ via
\[
    \updelta(u,v) \,\defeq\, \sum_{n=0}^\infty \updelta_n(u,v) \quad \text{for all } u \sim_G v.
\]
Observe that, by \eqref{eq:edge_bound1} and \eqref{eq:edge_bound}, $\|\updelta(u,v)\|_\infty \leq 1$ for all $uv \in E(G)$, and hence
\[
    \|\updelta(u,v)\|_\infty \,\leq\, \dist_G(u,v) \quad \text{for all } u \sim_G v.
\]

On the other hand, suppose that $u$, $v \in V(G)$ are two vertices such that $\dist_G(u,v) \in \interval_n$ for some $n \in \N$. Then there are at least $3m/8$ indices $i \in [m]$ satisfying
\begin{equation}\label{eq:partialsum}
    \left|\sum_{k = 0}^n \updelta_{k,i}(u,v)\right| \,=\, \left|\phi_{n+1,i}(v) - \phi_{n+1, i} (u)\right| \,\geq\, \dist_G(u,v)^{1-\epsilon}.
\end{equation}
Since $\D_{n+1}$ is a $3m/4$-strong $(r_{2n+2}, \alpha)$-padded decomposition, there are at least $3m/4$ indices $i \in [m]$ such that $B_G(u, r_{2n+2})$ is contained in a $\mP_{n+1,i}$-cluster. As $3/4 + 3/8 > 1$, there must exist an index $i$ such that $B_G(u, r_{2n+2})$ is contained in a $\mP_{n+1,i}$-cluster and \eqref{eq:partialsum} holds. Since
\[\dist_G(u,v) \,\leq\, R_n^\gamma \,=\, r_{2n+1}^{\alpha \gamma/\beta} \,=\, r_{2n+1}^{\alpha^2}  \,\leq\, r_{2n+1}^\beta \,=\, r_{2n+2},\]
with this choice of $i$, $u$ and $v$ must belong to the same $\mP_{n+1,i}$-cluster, and hence $\updelta_{k, i} (u,v) = 0$ for all $k > n$. Thus, we have $\sum_{k = 0}^\infty \updelta_{k,i}(u,v) = \sum_{k = 0}^n \updelta_{k,i}(u,v)$, and so $\norm{\updelta(u,v)}_\infty \geq \dist_G(u,v)^{1-\epsilon}$ by \eqref{eq:partialsum}. To summarize, we have shown that the inequality $\norm{\updelta(u,v)}_\infty \geq \dist_G(u,v)^{1-\epsilon}$ holds whenever
\[
    \dist_G(u,v) \,\in\, \bigcup_{n = 0}^\infty \interval_n \,=\, \bigcup_{n=0}^\infty \left(x^{\alpha \beta^{2n+1}}, \, x^{\alpha^2 \beta^{2n+1}}\right].
\]

To complete the proof, we let $\ell \in \N^+$ be the smallest integer such that $\alpha^\ell \geq \beta^2$, i.e.,
\[
    \ell \,\defeq\, \left\lceil \frac{2 \log \beta}{ \log \alpha}\right\rceil \,\leq\, \frac{500 \log(1/\epsilon)}{\epsilon}.
\]
Run the above construction for the following $\ell$ values of $x$:
\[
    x_0 \,\defeq\, \max\left\{\left(\frac{10^7 \, \bound}{\epsilon}\right)^{\frac{200}{\epsilon}}, \, r \right\}, \quad x_1 \,\defeq\, x_0^{\alpha}, \quad x_2 \,\defeq\, x_1^\alpha, \quad \ldots, \quad x_{\ell - 1} \,\defeq\, x_{\ell - 2}^\alpha.
\]
Let the resulting Borel $\Z^m$-cocycles be $\updelta^0$, \ldots, $\updelta^{\ell-1} \colon {\sim_G} \to \Z^m$. Define
\[
    \updelta \,\defeq\, \bigoplus_{j = 0}^{\ell - 1} \updelta^j \,\colon\, {\sim_G} \to \Z^{\ell m}.
\]
Then $\updelta$ is a Borel $\Z^{\ell m}$-cocycle such that for all $u \sim_G v$,
\[
    \|\updelta(u,v)\|_\infty \,=\, \max \set{\|\updelta^j(u,v)\|_\infty \,:\, 0 \leq j \leq \ell-1} \,\leq\, \dist_G(u,v).
\]
Furthermore, we have $\norm{\updelta(u,v)}_\infty \geq \dist_G(u,v)^{1-\epsilon}$ whenever
\begin{align*}
    \dist_G(u,v) \,&\in\, \bigcup_{n = 0}^\infty  \bigcup_{j = 0}^{\ell-1} \left(x_j^{\alpha \beta^{2n+1}}, \, x_j^{\alpha^2 \beta^{2n+1}}\right] \\
    \left[\text{since } x_{\ell-1} = x_0^{\alpha^{\ell-1}} \geq x_0^{\beta^2/\alpha}\right] \qquad &\supseteq\, \bigcup_{n = 0}^\infty  \left(\, \bigcup_{j = 0}^{\ell - 2}  \left(x_j^{\alpha \beta^{2n+1}}, \, x_{j+1}^{\alpha \beta^{2n+1}}\right] \,\cup\, \left(x_{\ell-1}^{\alpha \beta^{2n+1}}, \, x_{0}^{\alpha \beta^{2n+3}}\right] \right) \\
    &=\, \left(x_0^{\gamma}, \, +\infty\right).
\end{align*}
It follows that $\norm{\updelta(u,v)}_\infty \geq \dist_G(u,v)^{1-\epsilon}$ for all $u \sim_G v$ with
\[
    \dist_G(u,v) \, \geq\, \max \left\{\left(\frac{10^7 \, \bound}{\epsilon}\right)^{\frac{3000}{\epsilon^2}}, \, r^{\frac{15}{\epsilon}} \right\} \,=\, x_0^{\frac{15}{\epsilon}} \,>\, x_0^{\gamma}.
\]
Also, for every vertex $u \in V(G)$, all vertices $v \sim_G u$ such that $\updelta(u,v) = 0$ are contained in $B_G(u, x_0^{\gamma})$. Hence, the equivalence relation $E_\updelta$ is finite. Since finite Borel equivalence relations are smooth \cite[Lemma 5.21]{pikhurko2020borel}, the cocycle $\updelta$ is realizable by Theorem \ref{thm:cocycle-embedding}. Therefore, by Lemma \ref{lemma:map_to_free}, there exists a Borel map $f \colon V(G) \to \Free(2^{\Z^{\ell m + 1}})$ such that for all $u$, $v \in V(G)$,
\[
    \dist_\infty(f(u), f(v)) \,=\, \norm{\updelta(u,v)}_\infty.
\]
It remains to observe that
\[
    \ell m + 1 \,\leq\, \frac{500 \log(1/\epsilon)}{\epsilon} \, \left\lceil \frac{1440 \bound}{\epsilon}\right\rceil \,+\, 1 \, \leq\, \frac{10^7 \, \bound \log(1/\epsilon)}{\epsilon^2},
\]
and hence the function $f$ is as desired.
\end{proof}

\section{Injective coarse embeddings}\label{sec:injective}

    In this section we use Theorem~\ref{theo:coarse_p} to construct injective coarse embeddings of graphs of polynomial growth into grids with appropriately chosen parameters: 


  \begin{theo}\label{theo:injcoarse1_p}
        Let $G$ be a Borel $(\bound, r)$-graph with $\bound$, $r \geq 1$. Let $s \in \N$ and $ 0 < \epsilon < 1/2$ and define
        \[
            N \,\defeq\, \left\lceil \frac{10^7 \, \bound \log(1/\epsilon)}{\epsilon^2} \right\rceil,  \quad R \,\defeq\, \max \left \{ \left(\frac{10^7 \, \bound}{\epsilon}\right)^{\frac{3000}{\epsilon^2}}, \, r^{\frac{15}{\epsilon}}\right\}, \quad \text{and} \quad k \,\defeq\, \left\lceil\frac{\bound \log R}{\log(s+1)}\right\rceil.
        \]
        Then 
        there exists an injective Borel function $f \colon V(G) \to \Free(2^{\Z^{N+k+1}})$ 
        such that for all $u$, $v \in V(G)$, 
        \begin{itemize}
            \item $\dist_\infty(f(u), f(v)) \leq \max \set{\dist_G(u,v), s}$, and
            \item if $\dist_G(u,v) \geq R$, 
        then $\dist_\infty(f(u), f(v)) \geq \dist_G(u,v)^{1-\epsilon}$.
        \end{itemize}
        %
    \end{theo}
\begin{proof}
    By Theorem~\ref{theo:coarse_p}, there is a Borel function $g \colon V(G) \to \Free(2^{\Z^N})$ such that for all $u$, $v \in V(G)$, we have $\dist_\infty(g(u), g(v)) \leq \dist_G(u,v)$ and if $\dist_G(u,v) \geq R$, then $\dist_\infty(g(u), g(v)) \geq \dist_G(u,v)^{1-\epsilon}$. The latter condition implies that $g$ is a Borel reduction from $\sim_G$ to the orbit equivalence relation of the shift action $\Z^N \acts \Free(2^{\Z^N})$.
    We let $\updelta g \colon {\sim_G} \to \Z^N$ denote the cocycle associated to $g$. 
    
    Note that if $g(u) = g(v)$, then $\dist_G(u,v) <R$. Thus, our goal is to use $k + 1$ additional coordinates to distinguish the vertices at distance less than $R$ from each other. To this end, let $G^R$ be the graph with vertex set $V(G)$ in which two vertices $u$, $v$ are adjacent if and only if $1 \leq \dist_G(u,v) \leq R$. The maximum degree of $G^R$ satisfies \[\Delta(G^R) \leq \gamma_G(R) - 1 \leq R^\bound - 1.\] Therefore, by Theorem~\ref{theo:Delta+1}, the Borel chromatic number of $G^R$ is at most $R^\bound$. Note that $R^\bound \leq (s+1)^k$ by the choice of $k$, so we can fix a Borel proper coloring
        \[
            c \colon V(G) \to \set{0,1,\ldots, s}^k \subset \Z^k
        \]
        of $G^R$. In other words, $c$ is a Borel function such that $c(u) \neq c(v)$ whenever $1 \leq \dist_G(u,v) \leq R$. Now we define a $\Z^{N+k}$-cocycle $\updelta$ on $\sim_G$ as follows:
        \[
            \updelta(u,v) \,\defeq\, (\updelta g (u,v), \, c(v) - c(u)).
        \]
        Observe that for all $u$, $v \in V(G)$,
        \[
            \|\updelta(u,v)\|_\infty \,=\, \max\left\{\|\updelta g(u,v) \|_\infty, \, \|c(v) - c(u)\|_\infty\right\} \,=\, \max \left\{\dist_\infty (g(u), g(v)), \, \|c(v) - c(u)\|_\infty\right\}.
        \]
        Since $c(u)$, $c(v) \in \set{0,1,\ldots, s}^k$ and $g$ is a contraction, it follows that
        \[
            \|\updelta(u,v)\|_\infty \,\leq\, \max \{ \dist_G(u,v), s\}.
        \]
        If $1 \leq \dist_G(u,v) \leq R$, then $c(u) \neq c(v)$, and thus $\updelta(u,v) \neq 0$. And if $\dist_G(u,v) \geq R$, then
        \[
            \|\updelta(u,v) \|_\infty \,\geq\, \dist_\infty(g(u), g(v)) \,\geq\, \dist_G(u,v)^{1-\epsilon}.
        \]
        In particular, $\updelta(u,v) = 0$ if and only if $u = v$. In other words, the equivalence relation $E_\updelta$ is trivial, hence smooth, and thus the cocycle $\updelta$ is realizable by Theorem \ref{thm:cocycle-embedding}. Therefore, by Lemma \ref{lemma:map_to_free}, there exists a Borel map $f \colon V(G) \to \Free(2^{\Z^{N + k + 1}})$ such that for all $u$, $v \in v(G)$, \[\dist_\infty(f(u), f(v)) \, = \,\norm{\updelta(u,v)}_\infty.\] This function $f$ clearly has all the desired properties.
\end{proof}

    Theorems~\ref{theo:injcoarse1} and \ref{theo:injcoarse2} now follow from Theorem~\ref{theo:injcoarse1_p} by choosing specific values for the parameter $s$, namely $s = \ceil{R}$ for Theorem~\ref{theo:injcoarse1} and $s = 1$ for Theorem~\ref{theo:injcoarse2}.

\phantomsection
\addcontentsline{toc}{section}{References}

\printbibliography 
\end{document}